\documentclass[article,ij4uq]{ij4uq}              %

\usepackage[scaled]{helvet}

\usepackage[hang]{footmisc}
\usepackage{tikz}
\usepackage{tabu}
\usepackage{booktabs}%
\usetikzlibrary{calc}
\newcommand*{\myalign}[2]{\multicolumn{1}{#1}{#2}}
\fancypagestyle{plain}{%
  \fancyhf{}

  \fancyfoot[R]{\small\bf\thepage }
  \fancyfoot[L]{\fottitle}
  }
\renewcommand{\dmyy}{17}
\renewcommand{\myyear}{2017}
\renewcommand{\today}{}

\newcommand{\normm}[2][]{\left\Vert #2\right\Vert_{#1}}
\newcommand{\inprod}[2][]{\left\langle #2\right\rangle_{#1}}
\newcommand{\eps}{\varepsilon}
\newcommand{\bs}[1]{\ensuremath{\boldsymbol{#1}}}

\def\dontshow#1{{}}

\definecolor{grassgreen}{RGB}{92,135,39}

\begin{document}

\volume{Volume x, Issue x, \myyear\today}
\title{Estimation of the Robin coefficient field in a Poisson problem with uncertain conductivity field
}
\titlehead{}
\authorhead{R.Nicholson, N. Petra, \& J. P. Kaipio}
\corrauthor[1]{Ruanui Nicholson}
\author[2]{No\'{e}mi Petra}
\author[3]{Jari P. Kaipio}
\corremail{ruanui.nicholson@auckland.ac.nz}
\corraddress{Department of Engineering Science, University of Auckland, Private Bag 92019, Auckland Mail Center, Auckland 1142, New Zealand}
\address[1]{Department of Engineering Science, University of Auckland, Private Bag 92019, Auckland Mail Center, Auckland 1142, New Zealand}
\address[2]{School of Natural Sciences, University of California, Merced, 5200 North Lake Road, Merced, CA, 95343, United States}
\address[3]{Department of Mathematics, University of Auckland, Private Bag 92019, Auckland Mail Center, Auckland 1142, New Zealand}

\dataO{mm/dd/yyyy}
\dataF{mm/dd/yyyy}

\abstract{We consider the reconstruction of a heterogeneous
  coefficient field in a Robin boundary condition on an inaccessible
  part of the boundary in a Poisson problem with an uncertain (or
  unknown) inhomogeneous conductivity field in the interior of the
  domain.  To account for model errors that stem from the uncertainty
  in the conductivity coefficient, we treat the unknown conductivity
  as a nuisance parameter and carry out approximative
  premarginalization over it, and invert for the Robin coefficient
  field only. We approximate the related modelling errors via the
  Bayesian approximation error (BAE) approach.  The uncertainty
  analysis presented here relies on a local linearization of the
  parameter-to-observable map at the maximum a posteriori (MAP)
  estimates, which leads to a normal (Gaussian) approximation of the
  parameter posterior density. To compute the MAP point we apply an
  inexact Newton conjugate gradient approach based on the adjoint
  methodology. The construction of the covariance is made tractable by
  invoking a low-rank approximation of the data misfit component of
  the Hessian. Two numerical experiments are considered: one where the
  prior covariance on the conductivity is isotropic, and one where the
  prior covariance on the conductivity is anisotropic.  Results are
  compared to those based on standard error models, with particular
  emphasis on the feasibility of the posterior uncertainty
  estimates. We show that the BAE approach is a feasible one in the
  sense that the predicted posterior uncertainty is consistent with
  the actual estimation errors, while neglecting the related modelling
  error yields infeasible estimates for the Robin coefficient. In
  addition, we demonstrate that the BAE approach is approximately as
  computationally expensive (measured in the number of PDE solves) as
  the conventional error approach.}


\maketitle

\section{Introduction}\label{sect: intro}

In this paper, we consider the problem of estimating the unknown Robin
coefficient field in a Poisson problem with uncertain conductivity
field from available measurements.
This problem setup is inspired from physical applications, e.g., the
detection of corrosion of an electrostatic conductor \cite{Kaup1995,
  Inglese1997, Alessandrini2003} and the estimation of thermal
parameters \cite{Chantasiriwan1999,Divo2005}. The Poisson problem with
uncertain conductivity and Robin boundary coefficient has received considerable attention, both from
a theoretical standpoint as well as from a numerical point of view
\cite{Kaup1995, Inglese1997, Alessandrini2003, Chantasiriwan1999, Divo2005, Chaabane1999, Fasino1999, Chaabane2004, Jin2007, Ma2015}. However, the standard
assumption in these studies is that the internal conductivity (which
is generally a distribute parameter) is known. In this paper, we
consider both the conductivity and boundary condition to be uncertain. A
common approach would be to invert for both fields simultaneously,
however this results in a highly ill-posed and potentially
untracktable problem. To avoid the need for a joint inversion, we {\it
  premarginalize} over the internal conductivity, and then invert for
the Robin coefficient. Furthermore, to cope with the
infinite-dimensional nature of this inverse problem, we employ a
discretisation invariant method for the inversions
\cite{Bui-Thanh2013,Petra2014}. Thus the methods developed should be
immediately applicable to the case in which the Robin parameter is
high-dimensional (which is often the case for real applications).

There is a rich body of literature on theoretical and computational
aspects of the so-called inverse Robin problem, i.e. the problem of inferring the (distributed) Robin coefficient given measurements of the potential. In \cite{Inglese1997}, the
authors develop a direct reconstruction method based on a thin plate
approximation. In \cite{Chaabane1999}, several results on stability,
uniqueness and identifiability are established, while in
\cite{Alessandrini2003} a more general stability estimate is
proved. Numerical methods developed to solve the inverse Robin problem
include a quasi-reversibility method \cite{ Fasino1999} and an
approach based on an $L^1$-tracking functional
\cite{Chaabane2004}. More recently, in \cite{Jin2007}, a regularized
least-squares approach is taken via a variational formulation and in
\cite{Ma2015}, a regularised least-squares problem is solved using an
adjoint based approach, similar to the methods considered in the
present paper.

The previous studies consider the following inverse problem: given
noisy (partial point) measurements of $u$ on the boundary of a bounded
domain $\Omega$, with $\Omega\in\mathbb{R}^d$, $d\in\{2,3\}$,
determine the Robin coefficient field, $\beta(\bs{x})$. The field
$u$ satisfies the forward problem,
\begin{equation}\label{eq: laplace}
\begin{aligned}
-\Delta u({\bs x})&=0\quad&&\text{in }\Omega,\\
\nabla u({\bs x})\cdot {\bs n}_{\rm t}&=g({\bs x})\quad&&\text{on }\Gamma_{\rm t}\\
\nabla u({\bs x})\cdot {\bs n}_{\rm b}+\exp(\beta({\bs x})) u({\bs x})&=0\quad&&\text{on }\Gamma_{\!\rm b}\\
u({\bs x})&=0\quad&&\text{on }\Gamma_{\!\rm s},
\end{aligned}
\end{equation}
where (in the context of the present paper), $\Gamma_{\!\rm t}$ is
referred to as the top of the domain, $\Gamma_{\!\rm s}$ the sides of
the domain, and $\Gamma_{\!\rm b}$ the bottom of the domain. As such, we
have $\Gamma_{\rm t}\cap\Gamma_{\rm s}=\Gamma_{\rm s}\cap\Gamma_{\rm
  b}=\Gamma_{\rm t}\cap\Gamma_{\rm b}=\emptyset$, and (eventually
noisy pointwise) measurements of $u$ are available on $\Gamma_{\!\rm
  t}$. In the literature, Dirichlet boundary conditions are often
replaced by Neumann boundary conditions on $\Gamma_{\!\rm s}$, see for
example \cite{Chaabane2004,Jin2007}.
 
 However, in essentially all practical problems, the conductivity
 cannot be assumed to be a constant, leading to the spatially
 inhomogeneous problem
\begin{equation}\label{eq: simpleforward}
\begin{aligned}
-\nabla\cdot(\exp(a(\bs{x}))\nabla u({\bs x}))&=0\quad&&\text{in }\Omega,\\
\exp(a(\bs{x}))\nabla u({\bs x})\cdot {\bs n}_{\rm t}&=g({\bs x})\quad&&\text{on }\Gamma_{\!\rm t}\\
\exp(a(\bs{x}))\nabla u({\bs x})\cdot {\bs n}_{\rm b}+\exp(\beta({\bs x})) u({\bs x})&=0\quad&&\text{on }\Gamma_{\!\rm b}\\
u({\bs x})&=0\quad&&\text{on }\Gamma_{\!\rm s}.
\end{aligned}
\end{equation}
For example, in the case of detecting corrosion of an electrostatic
volume conductor, $\exp(a)$ models the electrical conductivity and it
is generally an unknown (distributed) parameter.

The estimation of the Robin coefficient $\beta(\bs{x})$ under the
forward model (\ref{eq: simpleforward}) would typically necessitate
the simultaneous estimation of the conductivity $a(\bs{x})$. Such a
task carries with it several challenges: The ill-posedness of the
problem would be increased significantly, and there is the potential
for issues around the identifiability of $\beta(\bs{x})$.  In this
paper, however, our main concern is the drastically increased
computational cost associated with estimating the parameter
$a(\bs{x})$ which is distributed over the entire volume $\Omega$
rather than estimating $\beta(\bs{x})$ only which is distributed on
$\Gamma_{\rm b}$ only.

The approach in this paper is based on the (initially)
infinite-dimensional formulation of the inverse problem, the adjoint
method for the computation of the related first and second order
derivative information, and the subsequent Bayesian approximation
error approach.  A brief review of these topics is given below.

The infinite dimensional approach to Bayesian inverse problems that
was developed in \cite{Stuart2010} is receiving considerable
attention, and has been successfully applied to several inverse
problems. The method is particularly suited to the case when the
parameter of interest is high-dimensional (stemming from the
discretization of the the unknown infinite-dimensional parameter
field), and it ensures convergence under discretization. An efficient
computational framework was developed in
\cite{Bui-Thanh2013,Petra2014}, based on an adjoint approach
\cite{Gunzburger2003, Hinze2009}, to implement the theoretical work
put forward in \cite{Stuart2010}, and was applied, for example, to
global seismic inversions in \cite{Bui-Thanh2013} and ice sheet flow
inverse problems in \cite{Petra2014}. The approach has also
successfully been applied to inverse acoustic obstacle scattering
problems in \cite{Bui-Thanh2014}. The infinite-dimensional Bayesian
setup has also been applied to optimal experimental design (OED) for
Bayesian nonlinear inverse problems governed by partial differential
equations~\cite{Alexanderian2016}. The goal of the OED problem was to
find an optimal placement of sensors (for measurements) so as to
minimize the uncertainty in the inferred parameter field.

The Bayesian approximation error (BAE) approach
\cite{Kaipio2005,Kaipio2007} was originally used as a means to take
into account the modelling errors induced by the use of reduced order
models.  The approach is based on approximate premarginalization over
modelling errors, which refers to a process similar to the
marginalization over additive errors to obtain the likelihood.
However, a particularl advantage of this method is the ability to
approximately premarginalize also over parameters which are not of
primary interest. In the context of electrical impedance tomography
(EIT), the BAE approach has been used to simultaneously premarginalize
over the unknown domain shape and the contact impedances of the
electrodes \cite{Nissinen2009}. Furthermore, in
\cite{Kolehmainen2011}, the approach was used to premarginalize over
the distributed scattering coefficient in diffuse optical tomography
(DOT) and, in \cite{Mozumder2016}, the method was used to
premarginalize over both the scattering and absorption coefficients in
the context of fluorescence diffuse optical tomography (fDOT). The BAE
method has also been applied to X-ray tomography to premarginalize
over distributed parameters outside a region of interest
\cite{Kaipio2013}.

In the context of premarginalizing over distributed parameters, the
BAE approach has thus far only been used to marginalize over unknowns
defined on spatial dimensions at most equal to that of the primary
parameter of interest.  I this paper, we show that the BAE approach is
also feasible for the premarginalization over a distributed parameter
in the entire domain when the parameter of primary interest is defined
only on (a subset of) the boundary of the domain.  We also show that
the infinite-dimensional framework for inverse problems posed in the
Bayesian setting is an effective method for solving the so-called
Robin inverse problem under an unknown (distributed) conductivity.

The paper is organized as follows.  In Section~2, we review the
infinite-dimensional framework for inverse problems, the adjoint
method for computation of derivative information, the computation of
the maximum a posteriori estimate and the approximate posterior
covariance.  In Section~3, we review the Bayesian approximation error
approach and, in Section~4, we formulate the problem of estimating the
Robin coefficient in the case of an unknown conductivity.  In
Section~5, we consider two numerical experiments: a conductivity with
spatially isotropic smooth covariance and one with an anisotropic
smooth covariance.  The results are compared to those based on
standard error models, with particular emphasis on the feasibility of
the posterior uncertainty estimates.

\section{Background on the Bayesian Approach to Inverse Problems in Infinite Dimensions}\label{sec: BA2IDIP}

In this section, we give a brief review of the formulation of Bayesian
inverse problems following~\cite{Bui-Thanh2013,Petra2014} to an extent
that is relevant to the present paper.  To this end, consider the
problem of finding ${\beta}(\bs{x})\in\mathcal{H}\subset
L^2(\Omega_\beta)$, from observed measurements ${\bs d}^{\rm
  obs}\in\mathbb{R}^q$, with ${\beta}$ and ${\bs d}^{\rm obs}$ linked
by
\begin{align}\label{eq: P1}
{\bs d}^{\rm obs}=\bs{f}_{{a}_*}( \beta)+\bs{e},
\end{align}
where $\bs{f}_{a_*}:\mathcal{H}\rightarrow\mathbb{R}^q$ is the
parameter-to-observable map, and $\bs{e}$ represents additive errors
in the measurements. The slightly unconventional notation
$\bs{f}_{a_*}$ used for the so-called {\it parameter-to-observable
  map} will be explained in Section~\ref{section: BAE}.  The fact that
${\beta}$ is by assumption infinite-dimensional, presents several
challenges.  First, there is no Lebesgue measure in infinite
dimensions and thus we cannot define the conventional notion of a
probability density function, and hence Bayes' formula must be
interpreted through the Radon-Nikodym derivative. Second, any prior
measure assigned to parameters must ensure well-posedness of the
inverse problem, that is, allow for the computation of the posterior.
Third, the discretization of the problem must be consistent with the
infinite-dimensional structure of the
problem~\cite{Stuart2010,Bui-Thanh2013,Petra2014}.

In this paper, we take the prior to be a Gaussian measure,
$\mu_\beta=\mathcal{N}(\beta_*,\mathcal{C}_{\beta})$ on
$L^2(\Omega_\beta)$, where $\beta_*$ is the prior mean, which lives in
$\mathcal{H}$, and $\mathcal{C}_{\beta}$ is the prior covariance
operator. As outlined in~\cite{Stuart2010}, the prior must be chosen
to satisfy certain regularity assumptions to ensure the Bayesian
inverse problem is well-defined. We employ a {\it weighted} squared
inverse elliptic operator as our prior covariance operator
\cite{Daon2016}, with the addition of homogeneous Robin (or Neumann)
boundary conditions. This is a slight modification to that used in
\cite{Bui-Thanh2013,Petra2014}, with the aim of mitigating any
artefacts in estimates or prior samples due to the
enforcement of boundary conditions. Specifically,
for $s\in$ $L^2(\Omega_\beta)$, the weak solution of
$\mathcal{A}\beta=s$ satisfies
\begin{align}\label{eq variationalprior}
\alpha_\beta\int_{\Omega_\beta}\bs{\gamma}_{\beta}\nabla \beta\cdot\nabla v+\beta v\;d\bs{y}+\int_{\partial\Omega_\beta}\kappa_\beta \beta v\;d\bs{t}=\int_{\Omega_\beta} s v\;d{\bs y}\quad\text{for all } v\in H^1(\Omega_\beta),
\end{align}
where $\alpha_\beta>0$ is inversely proportional to the prior variance, $\bs{\gamma}_{\beta}$ is a symmetric
positive definite uniformly bounded matrix controlling the correlation \cite{Bui-Thanh2013}, and
$\kappa_\beta\ge0$.
Then we take the prior covariance operator to be
\begin{align}
\mathcal{C}_\beta=\mathcal{W}\mathcal{A}^{-2}\mathcal{W}
\end{align}
where
\begin{align}\label{eq: normvar}
\mathcal{W}:=\frac{\sigma_\beta}{\sqrt{\mathcal{G}(\bs{x},\bs{x})}}\quad\text{and}\quad \sigma_\beta:=\frac{{\rm Ga}(\nu)}{(4\pi)^{d/2}\gamma_\beta^\nu\alpha_\beta^{{2}}}.
\end{align} 
Here $\mathcal{G}$ is the Greens function corresponding to
$\mathcal{A}^{2}$ in $\Omega$, ${\rm Ga}$ denotes the Gamma function,
$d$ is the spatial dimension, and $\nu+d/2=2$ \cite{Daon2016}.  The
addition of the weights normalizes the variance across the domain as
discussed below in Section~\ref{sec:numerics}. For efficient methods to extract
$\mathcal{G}(\bs{x},\bs{x})$, see for example, \cite{Rue2007,
  Bekas2007, Bekas2009, Lin2009, Tang2012}. We also note that there
are other methods aimed at mitigating the boundary effects, see for
example, \cite{Calvetti2006,Roininen2014,Daon2016}.

In this paper, we consider a normal noise model $\bs{e}\sim\mu_{\rm
  noise}=\mathcal{N}(0,{\bs \Gamma}_{\bs e})$, which results in the
likelihood
\begin{align}\label{eq: likelihood}
\pi_{\rm like}(\bs{d}^{\rm obs}|\beta)\propto\exp\left\{-\frac{1}{2}\left(\bs{f}(\beta)-\bs{d}^{\rm obs}\right)^T\Gamma_{\bs e}^{-1}\left(\bs{f}(\beta)-\bs{d}^{\rm obs}\right)\right\}.
\end{align}
In infinite dimensions, the Bayes' theorem states that the
Radon-Nikodym derivative of the posterior measure $\mu_{\rm d}$ with
respect to the prior measure $\mu_\beta$ is proportional to the
likelihood
\begin{align}\label{eq: radon_niko}
\frac{d\mu_{\rm d}}{d\mu_\beta}=\frac{1}{C}\pi_{\rm like}(\bs{d}^{\rm obs}|\beta)
\end{align}
where $C=\int_{\mathcal{H}}\pi_{\rm like}(\bs{d}^{\rm
  obs}|\beta)\,d\mu_\beta$ acts as a normalization
constant~\cite{Stuart2010}.
\subsection{Discretization of Bayesian Inverse Problems}
\label{sect: discret}
In this section, we review the finite-dimensional approximation of the
prior and the posterior distributions.  The discussion here follows
\cite{Petra2014, Bui-Thanh2013, Alexanderian2016}. The objective here
is to motivate the choice of the mass-weighted inner product space as
the correct space to work in, and to indicate the consequences.
Firstly, let $V_h$ denote a finite-dimensional subspace of
$L^2(\Omega_\beta)$ induced by a finite element discretization with
continuous Lagrange basis functions denoted by
$\{\phi_j\}_{j=1}^n$. The parameter of interest, $\beta\in
L^2(\Omega_\beta)$ is then approximated as $\beta_h=\sum_{j=1}^n
\beta_j\phi_j\in V_h$ and we then seek to invert for
$\bs{\beta}=[\beta_1,\beta_2,\dots,\beta_n]^T\in\mathbb{R}^n$.

We consider a Gaussian prior measure defined on $L^2(\Omega_\beta)$ and
thus the finite-dimensional subspace $V_h$ is equipped with the $L^2$
inner product. Thus any inner product between nodal coefficients will
be weighted by a mass matrix $\bs{M}$ so as to correctly approximate
the infinite-dimensional $L^2$ inner product. We denote the
mass-weighted inner product by $\inprod[{\bs M}]{\cdot,\cdot}$, with
$\inprod[{\bs M}]{\bs{y},\bs{z}}=\bs{y}^T{\bs M}\bs{z}$ and the 
symmetric positive definite mass matrix given by
\begin{align}
M_{ij}=\int_{\Omega_\beta}\phi_i(\bs{y})\phi_j(\bs{y})\;d\bs{y}\quad i,j\in\left\{1,2,\dots,n\right\}.
\end{align}
To distinguish between the Euclidean space $\mathbb{R}^n$ and the
$\mathbb{R}^n$ endowed with the mass-weighted inner product, we
introduce the notation $\mathbb{R}^n_{{\bs{M}}}$ to denote
$\mathbb{R}^n$ equipped with the mass-weighted inner product.
 
There are several crucial (yet subtle) differences which must be taken
into account when working in the $\mathbb{R}^n_{{\bs{M}}}$ rather than
in $\mathbb{R}^n$, which revolve largely around the concepts of
adjoints and transposes. For instance, for an operator
$\bs{B}:\mathbb{R}^n_{{\bs{M}}}\rightarrow \mathbb{R}^n_{{\bs{M}}}$,
we denote the matrix transpose by $\bs{B}^T$. The mass-weighted inner
product adjoint, however, denoted here by $\bs{B}^*$, satisfies
$\inprod[\bs{M}]{\bs{B}\bs{y},\bs{z}}=\inprod[\bs{M}]{\bs{y},\bs{B}^*\bs{z}}$,
for $\bs{y},\bs{z}\in\mathbb{R}^n_{\bs{M}}$, implying
\begin{align}
\bs{B}^*=\bs{M}^{-1}\bs{B}^T\bs{M}.
\end{align}
We also require two further adjoint type operators below in
Section~\ref{section: MapAppPost}. For this, let $\mathbb{R}^q$ and
$\mathbb{R}^r$ for some $q,r$, be endowed with the Euclidean inner
product.  We then define the adjoints $\bs{F}^\natural$ of
$\bs{F}:\mathbb{R}^n_{\bs{M}}\rightarrow\mathbb{R}^q$ and
$\bs{V}^\diamond$ of
$\bs{V}:\mathbb{R}^r\rightarrow\mathbb{R}^n_{\bs{M}}$ as
\begin{align}
\bs{F}^\natural&=\bs{M}^{-1}\bs{F}^T,\label{eq: transpose1}\\
\bs{V}^\diamond&=\bs{V}^T\bs{M}.
\end{align}

With these definitions in hand, the finite-dimensional approximation
of the operator $\mathcal{A}$ in (\ref{eq variationalprior}) is
$\bs{A}=\bs{M}^{-1}\bs{K}\bs{G}^{-1}$ where
\begin{align}
K_{ij}&=\alpha_\beta\int_{\Omega_\beta}\left(\theta_\beta\nabla\phi_i
\cdot\nabla\phi_j +\phi_i \phi_j
\right)\;d\bold{y}+\int_{\partial\Omega_\beta}\kappa_\beta\phi_i\phi_j\;d\bold{t},\\ G_{ij}&=\frac{1}{\sigma_\beta}\sqrt{K^{-1}_{ij}}\delta_{ij},\quad
i,j\in\{1,2,\dots,n\},
\end{align}
where $\sigma_\beta$ is as in (\ref{eq: normvar}) and $\delta_{ij}$ is the Kronecker delta.

We can now express the finite-dimensional approximation $\mu_\beta^h$
of the prior Gaussian measure, $\mu_\beta$, as
\begin{align}
\pi_{\rm prior}\left(\boldsymbol{\beta}\right)\propto\exp\left(-\frac{1}{2}\normm[\bold{M}]{\bold{A}\left(\boldsymbol{\beta}-\boldsymbol{\beta}_*\right)}^2\right),
\end{align}
where $\bs{\beta}_*$ is the discretization of the prior mean,
$\beta_*$, and the prior covariance matrix is by definition
$\bs{\Gamma}_{\bs{\beta}}=\bold{A}^{-2}$ (since $\bold{A}$ is
symmetric).  We can now state the finite-dimensional Bayes' formula
\begin{align}\label{eq: postproduct}
\pi_{\rm post}(\bs{\beta}):=\pi_{\rm post}(\bs{\beta}|\bs{d}^{\rm obs})\propto \pi_{\rm prior}(\bs{\beta})  \pi_{\rm like}(\bs{d}^{\rm obs}|\bs{\beta}),
\end{align}
where $\pi_{\rm post}(\bs{\beta})$ is the density of the
finite-dimensional approximation $\mu_d^h$ of the the posterior
measure $\mu_d$ and $\pi_{\rm like}$ is the likelihood given by
(\ref{eq: likelihood}), see \cite{Bui-Thanh2013,Petra2014}.  Thus we
can express the posterior density explicitly as
\begin{align}\label{eq: discpost}
\pi_{\rm post}(\bs{\beta})\propto \exp\left(  -\frac{1}{2}\normm[\bs{\Gamma}_{\rm e}^{-1}]{\bs{f}_{a_*}(\bs{\beta})-\bs{d}^{\rm obs} } ^2 -\frac{1}{2}\normm[\bold{M}]{\bold{A}\left(\boldsymbol{\beta}-\boldsymbol{\beta}_*\right)}^2 \right).
\end{align}

\subsection{The MAP estimate and the approximate posterior covariance}
\label{section: MapAppPost}
In the Bayesian framework, one seeks to determine the posterior
density of the parameter of interest.  In principle, then, one would
explore the posterior density typically with Markov chain Monte Carlo
sampling methods.  However, in large-scale problems with limited
computational resources, one often has to be content with a single
representative point estimate of the parameter along with an
(approximate) posterior covariance and credibility intervals.
Standard point estimates for the posterior include the conditional
mean (CM) and the maximum a priori (MAP) estimates.  In this paper, we
aim at computing the finite-dimensional MAP estimate $\bs{\beta}_{\rm
  MAP}$ and the Laplace (local normal) approximation for the
covariance which also yields approximate marginal distribution for
$\beta_{{\rm MAP},k}$.  For discussion on the extension of the MAP
estimate to infinite dimensions, see, for
example~\cite{Dashti2013,Helin2015}.

Following on from (\ref{eq: discpost}), the MAP estimate is defined as
the point in parameter space that maximizes the posterior probability
density function \cite{Tarantola2004,Kaipio2005}, that is,
\begin{align}\label{eq: map-min}
\bs{\beta}_{\rm MAP}:=\min_{\bs{\beta}\in \mathbb{R}^{n}}\frac{1}{2}\normm[\bs{\Gamma}_{\rm e}^{-1}]{\bs{f}_{a_*}(\bs{\beta})-\bs{d}^{\rm obs} } ^2 +\frac{1}{2}\normm[\bold{M}]{\bold{A}\left(\boldsymbol{\beta}-\boldsymbol{\beta}_*\right)}^2.
\end{align}
In line with \cite{Petra2014}, in (\ref{eq: map-min}) and the
remainder of the paper, we denote by $\bs{f}(\bs{\beta})$ the
parameter-to-observable map evaluated at the finite element function
corresponding to the parameter vector $\bs {\beta}$.

Assuming that the parameter-to-observable map is Fr\'{e}chet
differentiable, we can linearize (\ref{eq: P1}) around
$\bs{\beta}_{\rm MAP}$ and discretize to obtain the affine
approximation
\begin{align}
\bs{d}^{\rm obs} \approx \bs{f}_{a_*}(\bs{\beta}_{\rm MAP})+\bs{F}_{a_*}(\bs{\beta}_{\rm MAP})(\bs{\beta}-\bs{\beta}_{\rm MAP})+\bs{e},
\end{align}
where $\bs{F}_{a_*}(\bs{\beta}_{\rm MAP})$ is the Fr\'{e}chet
derivative of $\bs{f}_{a_*}(\beta)$ with respect to $\beta$ evaluated
at $\bs{\beta}_{\rm MAP}$.  The resulting approximation for the
discrete posterior measure $\mu_d^h$ (as well as in the
infinite-dimensional case), of $\beta$ conditional on $\bs{d}^{\rm
  obs}$ is then necessarily a Gaussian measure, given by
$\mathcal{N}(\bs{\beta}_{\rm MAP}, \bs{\Gamma}_{\rm post})$, with
covariance matrix
\begin{align}
\bs{\Gamma}_{\rm post}=\left(\bs{F}_{a_*}^\natural\bs{\Gamma}_{\bs e}^{-1}\bs{F}_{a_*}+\bs{\Gamma}_{\bs{\beta}}^{-1}\right)^{-1},
\end{align}
where $\bs{F}_{a_*}^\natural$ is the adjoint of $\bs{F}_{a_*}$, see
(\ref{eq: transpose1}).  We also note that the (approximate) posterior
covariance matrix is the inverse of the Gauss-Newton Hessian of the
negative log posterior (referred to simply as the Hessian from this
point onwards), denoted by $\bs{H}$, that is,
\begin{align}
\bs{\Gamma}_{\rm post}=\bs{H}^{-1}.
\end{align}

{\bf Low rank approximation for the approximate posterior covariance matrix.}
For an efficient calculation of the MAP estimate, and efficient action
of the approximate posterior covariance to vectors and action of the
square root of the Hessian on vectors (the latter is needed for
computing samples from the approximate posterior), we apply a low rank
approximation of the Hessian as detailed
in~\cite{FlathWilcoxAkcelikEtAl11,Bui-Thanh2013,Petra2014} and summartized below.

We note that the Hessian of the negative log posterior~\eqref{eq:
  discpost} can be split into the sum of a data misfit term,
$\bs{H}_{\rm mis}$, and the inverse of the prior covariance,
$\bs{\Gamma}_{\bs{\beta}}^{-1}$. By factorizing the prior covariance
as $\bs{\Gamma}_{\bs{\beta}}=\bs{L}\bs{L}^*$, we can rewrite the
Hessian as
\begin{align}\label{eq: HessianFact}
\bs{H}=\bs{H}_{\rm mis}+\bs{\Gamma}_{\bs{\beta}}^{-1}=\bs{H}_{\rm mis}+\bs{L}^{-*}\bs{L}^{-1}=\bs{L}^{-*}\left(\bs{L}^{*}\bs{H}_{\rm mis}\bs{L}+\bs{I}\right)\bs{L}^{-1}.
\end{align}
The final form of~\ref{eq: HessianFact} allows for an efficient method
of approximating $\bs{\Gamma}_{\rm post}$. The procedure relies on
exploiting the discretization invariant and often low rank nature of
$\bs{H}_{\rm
  mis}$~\cite{FlathWilcoxAkcelikEtAl11,Bui-Thanh2013,Petra2014,vogel2002}. Thus,
in this paper, we consider a low rank approximation for the so-called
{\it prior-preconditioned data misfit Hessian}~\cite{FlathWilcoxAkcelikEtAl11} which
takes the form
\begin{align}\label{eq: HessianFact2}
\bar{\bs{H}}_{\rm mis}=\bs{L}^{*}\bs{H}_{\rm mis}\bs{L}\approx\bs{V}_r\bs{\Lambda}_r\bs{V}^\diamond_r,
\end{align}
where $\bs{V}_r\in\mathbb{R}^{n\times r}$ contains the $r$ eigenvectors of the prior-preconditioned data misfit Hessian corresponding to the $r$ largest eigenvalues\footnote{The truncation value $r$ is chosen such that $\lambda_i\ll1$ for $i>r$  \cite{FlathWilcoxAkcelikEtAl11,Bui-Thanh2013,Petra2014}.} $\lambda_i$, $i=1,2,\dots,r$, and $\bs{\Lambda}_r=\text{diag}(\lambda_1,\lambda_2,\dots,\lambda_r)$.  
Then, by using the Sherman-Morrison-Woodbury formula \cite{Golub2013}, we have 
\begin{align}\label{eq: svdinvH}
\bs{H}^{-1}\approx\bs{L}\left(\bs{I}-\bs{V}_r\bs{D}_r\bs{V}^\diamond_r\right)\bs{L}^*,
\end{align}
where $\bs{D}_r=\text{diag}(\lambda_1/(\lambda_1+1),\lambda_2/(\lambda_2+1),\dots,\lambda_r/(\lambda_r+1))\in\mathbb{R}^{r\times r}$. Furthermore, for drawing samples from the Gaussian approximation for the posterior, we have 
\begin{align}\label{eq: SQRTsvdinvH}
\bs{H}^{-1/2}\approx\bs{L}\left(\bs{V}_r\bs{P}_r\bs{V}^\diamond_r+\bs{I}\right)\bs{M}^{-1/2},
\end{align}
where $\bs{P}_r=\text{diag}(\lambda_1/\sqrt{\lambda_1+1}-1,\lambda_2/\sqrt{\lambda_2+1}-1,\dots,\lambda_r/\sqrt{\lambda_r+1}-1)\in\mathbb{R}^{r\times r}$, see \cite{Bui-Thanh2013}.

\section{Background on Model discrepancy and the Bayesian approximation error approach}
\label{section: BAE}
Solving inverse problems in any framework, especially in the
large-scale case, can be computationally prohibitively expensive.  To
overcome this challenge, several classes of reduced-order and
surrogate models have been introduced. In \cite{Frangos2010}, methods
to reduce the computational cost of solving inverse problems in the
statistical setting are divided into three broad methods: reducing the
cost of a forward simulation, reducing the dimension of the input
parameter space, and reducing the number of samples required.  In line
with our goal of keeping the computational cost tracktable, in this
paper we focus on approximations of the posterior rather than on
sampling-based methods.

As an approach, the Bayesian approximation error approach (BAE),
developed in \cite{Kaipio2005,Kaipio2007}, can be seen to lie in the
first and second category. However, along with lowering computational
costs of forward simulations and reducing the dimension of the
parameter space, the BAE method can simultaneously take into account a
vast array of uncertainties in the forward model, see, for example,
\cite{Nissinen2011,Arridge2006,Kaipio2007a,Huttunen2007,Nissinen2008,Kolehmainen2011,Lipponen2011,Kaipio2013,Koponen2014,Mozumder2016}.
Below, we summarize the concept and implementation of the BAE method.

In the BAE, any errors induced by the use of surrogate models,
reduction of the parameter dimension, and/or model uncertainties are
propagated to a single additive error term.  Hence the form of the
posterior will be as in (\ref{eq: postproduct}) with a redefined
likelihood density. In what follows, let $\beta(\bs{x})$ be our
parameter of interest, and take $a(\bs{x})$ to denote a secondary
(nuisance) parameter.  By secondary, we mean that we do not wish to
estimate the unknown $a(\bs{x})$ but attempt to take the related
uncertainty into account and propagate the effects into the estimate
for the parameter of interest and the posterior uncertainty.  Except
for jointly normal linear models, it is not possible to {\it exactly
  premarginalize} over $a(\bs{x})$ \cite{Kaipio2005,Kaipio2013}.  In
the following, we outline how one can approximately premarginalize
over $a(\bs{x})$.  To this end, let
 \begin{align}
(a,\beta)\mapsto\bs{f}(a,\beta)
\end{align}
denote an accurate forward model, and let $\bs{e}$ again denote noise
which is additive and mutually independent with $\beta$ and $a$ such
that $\bs{e}\sim\mu_{\rm noise}=\mathcal{N}(0,{\bs \Gamma}_{\bs
  e})$. Then the accurate relationship between the parameters and
measurements is
 \begin{align}\label{eq: accurate}
{\bs d}^{\rm obs}&=\bs{f}(a,\beta)+{\bs e}.
\end{align}

In the BAE approach, rather than using the accurate model
$\bs{f}(a,\beta)$, we instead set $a=a_*$ and use the approximate
forward model
 \begin{align}
\beta\mapsto\bs{f}_{a_*}(\beta).
\end{align}
We note that in many applications the dimension of the parameter of interest is also reduced by projecting onto some reduced basis, see, for example, \cite{Kaipio2013} for more details. 
In general, replacing the accurate model with the approximate model introduces what has become known as {\it model discrepancy}, the difference between the predictions of the two models. To take into account this model discrepancy we rewrite  (\ref{eq: accurate}) as
 \begin{align}
 {\bs d}^{\rm obs}&=\bs{f}(a,\beta)+{\bs e}=\bs{f}_{a_*}(\beta)+{\bs e}+\underbrace{\left({\bs f}(a,\beta)-{\bs f}_{a_*}(\beta)\right)}_{=\boldsymbol{\eps}(a,\beta)}={\bs f}_{a_*}(\beta)+{\bs e}+\boldsymbol{\eps}=\bs{f}_{a_*}(\beta)+\boldsymbol{\nu},
\end{align}
where the discrepancy in the models, $\boldsymbol{\eps}$, 
is a random variable with the same dimensions as the measurements, and is known as the {\it approximation error} \cite{Kaipio2005, Kaipio2007,Nissinen2011}.
The sum $\boldsymbol{\nu}={\bs e}+\boldsymbol{\eps}$ is called the {\it total error} here.

At this point in the BAE approach, the following approximation is made:
\begin{align}
\bs{\eps}\vert \beta\sim\mu_{ \bs{\eps}|\beta} 
\approx \mathcal{N}(\bs{\eps}_{*\vert\beta},{\bs \Gamma}_{\bs \eps|\beta}),
\end{align}
that is, the conditional density of the approximation error $\bs{\nu}$
given the parameter of interest $\beta$ is approximated as
Gaussian\footnote{There is some work on retaining the full conditional
  density, see, for example, \cite{Calvetti2014,Calvetti2017}.}.  The
computation of $\bs{\eps}_{*\vert\beta}$ and ${\bs \Gamma}_{\bs
  \eps|\beta}$ is outlined in Section~\ref{sec:numerics}.  The fact
that the approximation error depends on $\beta$ implies that formally
$\boldsymbol{\eps}$ and $\beta$ cannot be taken as mutually
independent.  However, in several cases, such a further approximation
of independence leads to similar estimates for a significantly smaller
cost as explained in~\cite{Kaipio2013}.  With this further approximation, we have
\begin{align}
\boldsymbol{\nu}_*={\bs e}_*+\boldsymbol{\eps}_*\quad\text{and}\quad \bs{\Gamma}_{\nu}=\bs{\Gamma}_{e}+\bs{\Gamma}_{{\eps}},
\end{align}
which was originally referred to as the {\it enhanced error model} in \cite{Kaipio2005,Kaipio2007}%
\footnote{The actual form of $\bs{\eps}_{*\vert\beta}$ is
  $\bs{\eps}_{*\vert\beta} = \bs{\eps}_\ast +
  \bs\Gamma_{\eps\!\beta}\bs\Gamma_{\beta\!\beta}^{-1} (\beta -
  \beta_\ast)$ which incorporates the full covariance structure of the
  normal approximation for $\pi(\bs{\eps},\beta)$.  The prior
  covariance of $\beta$ cannot, however, be used in place of
  $\bs\Gamma_{\beta\!\beta}$ in this conditional expectation. Rather,
  it must be based on the same draws as those used to compute
  $\bs\Gamma_{\eps\!\eps}$ as in Section~5 below. In practise, this
  leads to a a semidefinite estimate for $\bs\Gamma_{\beta\!\beta}$,
  and the associated rank-deficient forms for the conditional
  expectations need to be employed \cite{Kaipio2013}.}.

The BAE approach results in both a revised functional, which the MAP
estimate minimizes, and a reformulated approximate posterior
covariance matrix. Specifically, we now have
\begin{align}
\bs{\beta}_{\rm MAP}&=\min_{\bs{\beta}\in \mathbb{R}^{n}}\frac{1}{2}\normm[\bs{\Gamma}_{ \nu}^{-1}]{\bs{f}_{a_*}(\bs{\beta})-\bs{d}^{\rm obs}+\bs{\nu}_* } ^2 +\frac{1}{2}\normm[\bold{M}]{\bold{A}\left(\boldsymbol{\beta}-\boldsymbol{\beta}_*\right)}^2,\label{bae map}\\
\bs{\Gamma}_{\rm post}&=\left(\bs{F}_{a_*}^\natural\bs{\Gamma}_{ \nu}^{-1}\bs{F}_{a_*}+\bs{\Gamma}_{\beta}^{-1}\right)^{-1}\label{bae post}.
\end{align}
We note that the infinite-dimensional counterparts to~\eqref{bae map}
and \eqref{bae post} can be formulated naturally.

As an indicator as to whether or not inclusion of the approximation
errors is appropriate the following rule of thumb can be adopted
\cite{Kaipio2013}: if
\begin{align}
\normm[]{\bs{e}_*}^2+\text{trace}(\bs{\Gamma}_e)<\normm[]{\bs{\eps}_*}^2+\text{trace}(\bs{\Gamma}_{\eps})
\end{align}
holds, then the approximation errors {\it dominate} the noise and neglecting the approximation errors will generally result in meaningless reconstructions, as demonstrated in Section \ref{sec:Results}. Moreover, if 
\begin{align}
\bs{e}_*^2(k)+\bs{\Gamma}_{e}^2(k,k)<\bs{\eps}_*(k)+\bs{\Gamma}_{\eps}(k,k)
\end{align}
for any $k$, then neglecting the approximation errors can still lead to meaningless results \cite{Kaipio2013}.

\section{Recovery of the Robin Coefficient}
\label{sect: RecRC}
In this section, we formulate the inverse Robin problem with a
spatially varying conductivity coefficient $a(\bs{x})$ which will
later be interchanged for a fixed conductivity $a_\ast$ In the chosen
geometry, we refer to the Robin coefficient as ``basal" since this
condition is posed only on the bottom part of a slab.  The
measurements are taken to be pointwise (noisy) potential measurements
on the top of the domain, while premarginalizing over $a(\bs{x})$. We
solve the inverse problem with Newton's method. Therefore, in what
follows, we formulate the forward problem and derive the corresponding
first and second order adjoint problems for the gradient and the
action of the Hessian to a vector needed by the optimization method.

\subsection{The Forward Problem}
As a model problem, let
$\Omega=\left[0,L\right]\times\left[0,L\right]\times\left[0,H\right]\in\mathbb{R}^3$
with $0<H \ll L<\infty$ denote the domain of the problem (a thin slab)
with boundary $\partial\Omega$. In our regime, a flux is prescribed on
$\Gamma_{\rm
  t}:=\left[0,L\right]\times\left[0,L\right]\times\left\{H\right\}$,
while the potential is measured at $q$ points on $\Gamma_{\rm t}$, see
Figure \ref{fig: probsetup} for a schematic representation.  A
homogenous Robin boundary condition is prescribed on
$\Omega_\beta=\Gamma_{\rm
  b}:=\left[0,L\right]\times\left[0,L\right]\times\left\{0\right\}$
while on the remainder of the boundary, $\Gamma_{\rm
  s}:=\partial\Omega\setminus\left(\Gamma_{\rm t}\cap\Gamma_{\rm
  b}\right)$, homogeneous Dirichlet boundary conditions are specified.
The conductivity coefficient $\exp(a(\bs{x}))$ is taken to be
spatially distributed random field in $\Omega$, while the Robin
coefficient is taken to be spatially varying random field on~$\Gamma_{\rm b}$.
To summarize, the forward problem reads
\begin{equation}\label{eq: Poisson}
\begin{aligned}
-\nabla\cdot\left(\exp(a({\bs x}))\nabla u({\bs x})\right)&=0\quad&\text{in }&\Omega\\
\exp(a({\bs x}))\nabla u({\bs x})\cdot {\bs n}_{\rm t}&=g({\bs x})\quad&\text{on }&\Gamma_{\rm t}\\
\exp(a({\bs x}))\nabla u({\bs x})\cdot {\bs n}_{\rm b}+\exp(\beta({\bs x})) u({\bs x})&=0\quad&\text{on }&\Gamma_{\rm b},\\
u({\bs x})&=0\quad&\text{on }&\Gamma_{\rm s}.
\end{aligned}
\end{equation}
where $u({\bs x})$ the potential, $g({\bs x})$ is the flux through
$\Gamma_{\rm t}$ with unit normal ${\bs n}_{\rm t}$, and $\Gamma_{\rm
  b}$ has unit normal ${\bs n}_{\rm b}$.
\dontshow{Noiseless data comprises of point measurements of the (forward) solution $u({\bs x})$ 
at points throughout $\Gamma_{\rm t}$. }
We employ the finite element method (FEM) for the numerical
approximation of the forward problem, with the standard Lagrange
piecewise linear nodal basis functions.

\subsection{Adjoint-Based Gauss-Newton Method for Solving the Inverse Problem}
\label{section: adjointopt}
We employ an inexact Newton-CG approach to solve the minimization problem (\ref{eq: map-min}) which requires both the gradient and the Hessian of the negative log prior and likelihood. To avoid calculations of forward sensitivities, which would require as many forward solves as the dimension of the parameter, we employ the adjoint approach \cite{Gunzburger2003, Hinze2009, Troeltzsch2010, Borzi2012} to compute the (infinite-dimensional) derivatives, which we show next. 

We denote the observation operator with $\mathcal{B}$ so that the
parameter-to-observable mapping can be written as $\bs{f}_{a_*} =
\mathcal{B}u$. Hence the infinite-dimensional counterpart of the
functional to be minimized in (\ref{eq: map-min}) can be rewritten as
\begin{align}\label{eq: map-min2}
\mathcal{J}(\beta)=\frac{1}{2}\normm[\bs{\Gamma}_{\bs \nu}^{-1}]{\mathcal{B}u-\bs{d}^{\rm obs}+\bs{\nu}_*}^2+\frac{1}{2}\normm[L^2(\Gamma_{\rm b})]{\mathcal{A}\mathcal{W}^{-1}(\beta-{\beta}_*)}^2, 
\end{align}
where $u(\bs{x})$ solves the forward problem (\ref{eq:
  Poisson}). Furthermore, let us define the space,
\begin{align}
\mathcal{V}:=\left\{v\in H^1(\Omega) : \left.v\right|_{\Gamma_{\rm s}}=0\right\}, 
\end{align}
then we can define the Lagrangian functional $\mathcal{L}:\mathcal{V}\times \mathcal{V}\times \mathcal{E}\rightarrow\mathbb{R}$,
\begin{align}
\mathcal{L}(u,p,\beta):=\mathcal{J}(\beta)+\int_{\Omega}\exp(a(\bs{x}))\nabla u\cdot\nabla p\;d\bs{x} -\int_{\Gamma_{\rm t}}gp\;d\bs{s}_{\rm t}+\int_{\Gamma_{\rm b}}\exp(\beta) up\;d\bs{s}_{\rm b},
\end{align}
The space $\mathcal{E}$ is the Cameron-Martin space
$\mathcal{E}=\text{range}(\mathcal{C}^{\frac{1}{2}})=\text{dom}(\mathcal{A})$,
induced by the prior measure, see \cite{Stuart2010} for full details,
or for example \cite{Alexanderian2016} for a brief overview.

Determining the gradient of $\mathcal{J}$ is achieved by requiring that variations
of the Lagrangian $\mathcal{L}$ with respect to the forward potential $u$ and the so-called {\it adjoint potential} $p$ vanish.  
This results in the following 
strong form of the gradient $\mathcal{G}$ for the variations with respect to $\beta$
\begin{align}\label{eq: infdimGrad}
\mathcal{G}(\beta):=\mathcal{W}^{-1}\mathcal{A}^2\mathcal{W}^{-1}\left(\beta-\beta_*\right)+\exp(\beta)up,
\end{align}
with $u$ being the solution of the {\it forward Poisson problem} (\ref{eq: Poisson}) for given $\beta$, 
while $p$ satisfies the following {\it adjoint Poisson problem} for given $u(\bs{x})$ and $\beta(\bs{x})$
\begin{equation}\label{eq: AdjPoisson}
\begin{aligned}
-\nabla\cdot\left(\exp(a({\bs x}))\nabla p({\bs x})\right)&=-\mathcal{B}^*\bs{\Gamma}_{\bs{\nu}}^{-1}(\mathcal{B}u(\bs{x})-\bs{d}^{\rm obs}+\bs{\nu}_*)\quad&\text{in }&\Omega,\\
\exp(a({\bs x}))\nabla p({\bs x})\cdot {\bs n}_{\rm t}&=0\quad&\text{on }&\Gamma_{\rm t},\\
\exp(a({\bs x}))\nabla p({\bs x})\cdot {\bs n}_{\rm b}+\exp(\beta({\bs x})) p({\bs x})&=0\quad&\text{on }&\Gamma_{\rm b},\\
p({\bs x})&=0\quad&\text{on }&\Gamma_{\rm s}.
\end{aligned}
\end{equation}
The action of the Gauss-Newton approximation of the Hessian operator evaluated at $\beta$ in
the direction $\hat{\beta}$ is given by
\begin{align}\label{eq: infdimHess}
\mathcal{H}(\beta)(\hat{\beta}):=\mathcal{W}^{-1}\mathcal{A}^2\mathcal{W}^{-1}\hat{\beta}+\exp(\beta)\hat{\beta}u\hat{p},
\end{align}
where the {\it incremental adjoint potential} $\hat{p}$ satisfies the
so-called {\it incremental adjoint (or second order adjoint) Poisson
  problem}
\begin{equation}\label{eq: IncAdjPoisson}
\begin{aligned}
-\nabla\cdot\left(\exp(a({\bs x}))\nabla \hat{p}({\bs x})\right)&=-\mathcal{B}^*\bs{\Gamma}_{\bs{\nu}}^{-1}\mathcal{B}\hat{u}(\bs{x})\quad&\text{in }&\Omega\\
\exp(a({\bs x}))\nabla \hat{p}({\bs x})\cdot {\bs n}_{\rm t}&=0\quad&\text{on }&\Gamma_{\rm t}\\
\exp(a({\bs x}))\nabla \hat{p}({\bs x})\cdot {\bs n}_{\rm b}+\exp(\beta({\bs x})) \hat{p}({\bs x})&=0\quad&\text{on }&\Gamma_{\rm b},\\
\hat{p}({\bs x})&=0\quad&\text{on }&\Gamma_{\rm s},
\end{aligned}
\end{equation}
and the {\it incremental forward potential} $\hat{u}$ satisfies the
{\it incremental forward Poisson problem}
\begin{equation}\label{eq: IncPoisson}
\begin{aligned}
-\nabla\cdot\left(\exp(a({\bs x}))\nabla \hat{u}({\bs x})\right)&=0\quad&\text{in }&\Omega\\
\exp(a({\bs x}))\nabla \hat{u}({\bs x})\cdot {\bs n}_{\rm t}&=0\quad&\text{on }&\Gamma_{\rm t}\\
\exp(a({\bs x}))\nabla \hat{u}({\bs x})\cdot {\bs n}_{\rm b}+\exp(\beta({\bs x})) \hat{u}({\bs x})&=-\hat{\beta}\exp(\beta({\bs x})) u({\bs x})\quad&\text{on }&\Gamma_{\rm b},\\
\hat{u}({\bs x})&=0\quad&\text{on }&\Gamma_{\rm s}.
\end{aligned}
\end{equation}
The resulting system to be solved for the Gauss-Newton search
direction, $\hat{\beta}$, is
\begin{align}\label{eq: infdimGNdir}
\mathcal{H}(\beta)(\hat{\beta})=-\mathcal{G}(\beta).
\end{align}
In the case of a high-dimensional parameter, it may not be feasible to
solve the system (\ref{eq: infdimGNdir}) directly, and iterative methods
are usually employed.  To this end, we employ the conjugate gradient
(CG) method, which requires only the action of the (discrete) Hessian,
that is, Hessian-vector products. By examining (\ref{eq: infdimHess}),
we see that each Hessian-vector product amounts to solving both the
adjoint Poisson problem (\ref{eq: IncPoisson}) and the incremental
adjoint Poisson problem (\ref{eq: IncAdjPoisson}).  Moreover, the
computation of the gradient (\ref{eq: infdimGrad}) requires the
solutions of the forward Poisson problem (\ref{eq: Poisson}) and the
adjoint Poisson problem (\ref{eq: AdjPoisson}).  The dominant cost of
solving (\ref{eq: infdimGNdir})
is in computing the solution to Poisson equations, and thus the
computational cost of the method can be roughly measured in number of
Poisson problems solved.

Carrying out the premarginalization over $a$ avoids solving for $a$
since we replace $a(\bs{x})$ with $a_*(\bs{x})$.  Indeed, solving for
$a$ would require the computation of another search direction, which
in the Gauss-Newton case would be computed by solving
\begin{align}\label{eq: NewtonFor a}
\mathcal{H}_a(a)(\hat{a})=-\mathcal{G}_a(a)
\end{align}
for $\hat{a}$, with 
\begin{align}
\mathcal{G}_a(a):=\mathcal{A}_a^2\left(a-a_*\right)+\exp(a)\nabla u\cdot \nabla p\quad\text{and}\quad\mathcal{H}_a(a)(\hat{a}):=\mathcal{A}_a^2\hat{a}+\exp(a)\hat{a}\nabla u\cdot \nabla \hat{p}
\end{align}
the gradient and Hessian respectively, and $\mathcal{A}_a$ is defined
in (\ref{eq: priorona}).  In addition, the system (\ref{eq: NewtonFor
  a}) would need to be solved over the entire domain $\Omega$,
necessitating an extra order of magnitude of computational complexity.
Such an approach is outlined for the nonlinear Stokes flow problem in
\cite{Petra2012}.

\section{Numerical Examples}
\label{sec:numerics}
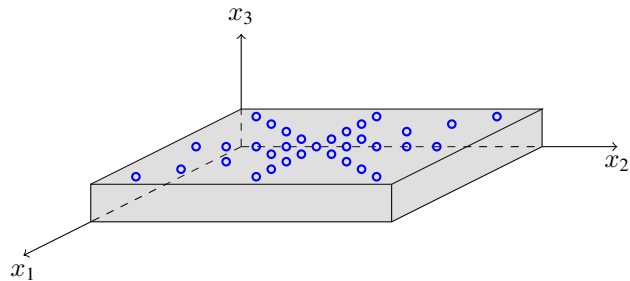
\begin{figure}[b!]
\centering
\begin{tikzpicture}

\fill[gray!99,nearly transparent] (0,0) -- (4,0) -- (4,.5) -- (0,.5) -- cycle;%
\fill[gray!99,nearly transparent] (4,0) -- (4,.5) -- (6,1.5) -- (6,1); -- cycle;%
\draw (0,0) -- (4,0) -- (4,.5) -- (0,.5) -- (0,0);%
\draw  (4,.5) -- (6,1.5) -- (2,1.5) -- (0,.5);%
\draw  (6,1.5) -- (6,1) -- (4,0);%
\draw[->] (2,1.5) -- (2,2.5) node[above] {${ x_3}$};
\draw[dashed] (2,1) -- (2,1.5) ;

\draw[->] (6,1) -- (7,1) node[below] {${ x_2}$};%
\draw[dashed] (2,1) -- (6,1);%

\draw[dashed] (2,1) -- (0,0);%
\draw[->] (0,0) -- (-2/2.23606,-1/2.23606) node[below] {${ x_1}$};%

\fill[gray!99,nearly transparent] (4,.5) -- (6,1.5) -- (2,1.5) -- (0,.5) -- cycle;
\draw [thick,blue] (2.2,2.5-4.4/4) circle [radius=0.05];
\draw [thick,blue] (2.4,2.5-4.8/4) circle [radius=0.05];
\draw [thick,blue] (2.6,2.5-5.2/4) circle [radius=0.05];
\draw [thick,blue] (2.8,2.5-5.6/4) circle [radius=0.05];
\draw [thick,blue] (3,2.5-6/4) circle [radius=0.05];
\draw [thick,blue] (3.2,2.5-6.4/4) circle [radius=0.05];
\draw [thick,blue] (3.4,2.5-6.8/4) circle [radius=0.05];
\draw [thick,blue] (3.6,2.5-7.2/4) circle [radius=0.05];
\draw [thick,blue] (3.8,2.5-7.6/4) circle [radius=0.05];

\draw [thick,blue] (.6,.5+.6/6) circle [radius=0.05];
\draw [thick,blue] (1.2,.5+1.2/6) circle [radius=0.05];
\draw [thick,blue] (1.8,.5+1.8/6) circle [radius=0.05];
\draw [thick,blue] (2.4,.5+2.4/6) circle [radius=0.05];
\draw [thick,blue] (3.6,.5+3.6/6) circle [radius=0.05];
\draw [thick,blue] (4.2,.5+4.2/6) circle [radius=0.05];
\draw [thick,blue] (4.8,.5+4.8/6) circle [radius=0.05];
\draw [thick,blue] (5.4,.5+5.4/6) circle [radius=0.05];

\draw [thick,blue] (2.2,-.5+2.2/2) circle [radius=0.05];
\draw [thick,blue] (2.4,-.5+2.4/2) circle [radius=0.05];
\draw [thick,blue] (2.6,-.5+2.6/2) circle [radius=0.05];
\draw [thick,blue](2.8,-.5+2.8/2) circle [radius=0.05];
\draw [thick,blue](3.2,-.5+3.2/2) circle [radius=0.05];
\draw [thick,blue](3.4,-.5+3.4/2) circle [radius=0.05];
\draw [thick,blue](3.6,-.5+3.6/2) circle [radius=0.05];
\draw [thick,blue](3.8,-.5+3.8/2) circle [radius=0.05];

\draw [thick,blue](1.4,1) circle [radius=0.05];
\draw [thick,blue](1.8,1) circle [radius=0.05];
\draw [thick,blue](2.2,1) circle [radius=0.05];
\draw [thick,blue](2.6,1) circle [radius=0.05];
\draw [thick,blue](3.4,1) circle [radius=0.05];
\draw [thick,blue](3.8,1) circle [radius=0.05];
\draw [thick,blue](4.2,1) circle [radius=0.05];
\draw [thick,blue](4.6,1) circle [radius=0.05];
\end{tikzpicture}
\caption{Set up for the model problem. A prescribed flux is set through the top of the domain and measurements of the potential are taken at points on the top of the domain (blue circles). A  
Robin boundary condition is enforced at the bottom of the domain, while the sides of the domain are prescribed homogeneous Dirichlet boundary conditions.}
\label{fig: probsetup}
\end{figure}
In this section, we consider two numerical experiments, one with a
conductivity with isotropic homogeneous covariance structure and one
with an anisotropic structure.  The latter structure is akin to
horizontally layered (stochastic) strata in which the correlation
length is smaller in the vertical direction than in the horizontal
plane. We will pay particular attention to the feasibility of the
posterior error estimates, that is, we will investigate whether the
posterior models (essentially) support the actual Robin coefficient.

\subsection{Problem setup}
In both experiments, the domain $\Omega\in\mathbb{R}^3$ is a
rectangular parallelepiped with thickness $H=0.01$ and width $L=1$,
such that $L/H=100$.  The measurements consists of $q=33$ point
measurements on the top of the domain, as illustrated in
Figure~\ref{fig: probsetup}.  To avoid the so-called inverse crime, we
use a finer FEM discretization to generate the synthetic data than the
FEM discretization used in the inversions. Moreover, the mesh used in
the second example to generate the data is finer than the
corresponding mesh used in the first example to ensure refinement of
the stratified conductivity in the volume. The details of the meshes
are presented in Table \ref{table: femMeshes}, in all cases Lagrange
piece-wise linear basis functions are used.

\begin{table}
\caption{Statistics of the mesh. The first column
  ({\bf Mesh used for}) relates what the mesh is used for; the second
  ({\fontseries{b}\selectfont\#}{\bf Nodes}) third
  ({\fontseries{b}\selectfont\#}{\bf Els}) and fourth
  ({\fontseries{b}\selectfont\#}{\bf Param}) columns give the number
  of total number of FEM nodes used in the entire volume, tetrahedral
  elements used in the entire volume, and the number nodes on the
  domain of the parameter, $\bs{\beta}$, respectively.}\label{table:
  femMeshes}
\noindent
\begin{tabu} to \textwidth {lXXXr}
   \toprule
   & {\bf Mesh used for} & {\fontseries{b}\selectfont\#}{\bf Nodes} &{\fontseries{b}\selectfont\#}{\bf Els}  &  {\fontseries{b}\selectfont\#}{\bf Param}\\
   \bottomrule
    Example 1 & & & &\\
  & \myalign{l}{Data synthesis} & \myalign{l}{28,611} & \myalign{l}{150,000} & \myalign{l}{2,601}\\
  & \myalign{l}{Inversion} & \myalign{l}{6,727} & \myalign{l}{32,400} & \myalign{l}{961}\\
     \bottomrule
    Example 2 & & & &\\
  & \myalign{l}{Data synthesis} & \myalign{l}{132,651} & \myalign{l}{750,000} & \myalign{l}{2,601}\\
  & \myalign{l}{Inversion} & \myalign{l}{29,791} & \myalign{l}{162,000} & \myalign{l}{961}\\
   \bottomrule
\end{tabu}
\end{table}

In both numerical examples, zero mean white noise is added to the
simulated measurements, with the noise covariance matrix given by
$\bs{\Gamma}_e=\delta^2_e\bs{I}$, with
$\delta_e=(\max(\mathcal{B}\bs{u})-\min(\mathcal{B}\bs{u}))\times1/100$,
that is, the noise level is $1\%$ of the range of the noiseless
measurements.

\begin{figure}[t!]
\includegraphics[width=\linewidth]{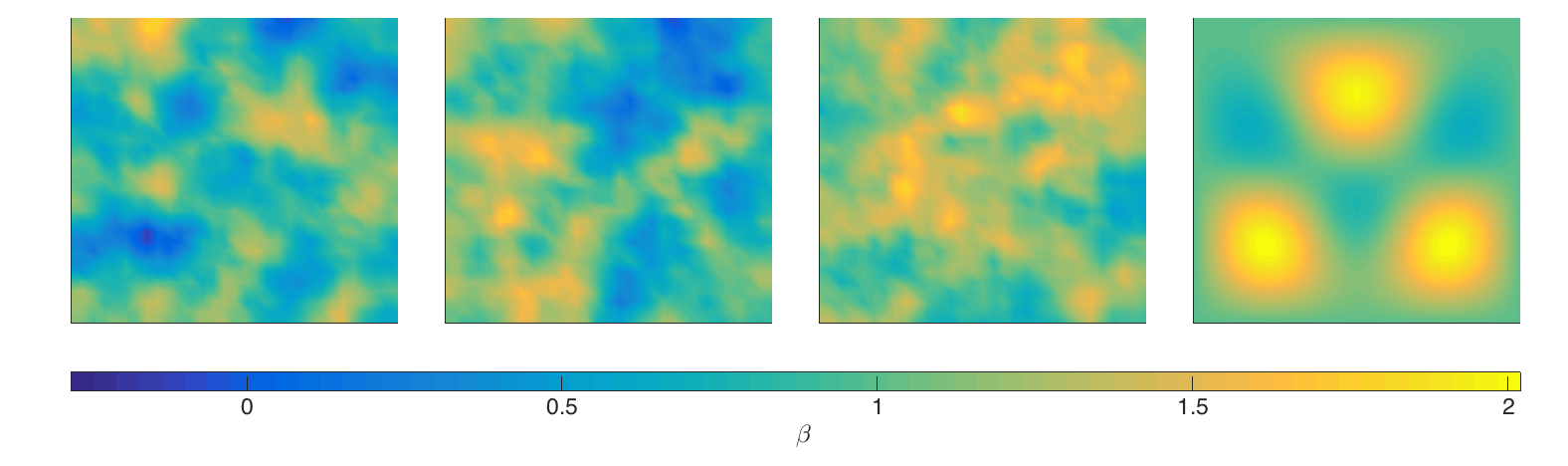}
\caption{Far left to centre right: Samples from the prior distribution on  
$\beta$. Right: The true value $\beta_{\rm true}$ used to compute the data.}
\label{fig: betadraws}
\end{figure}

{\bf Prior Models.} The prior density (normal random field) imposed on
$\beta$ is the same in both experiments, as outlined in Section
\ref{sec: BA2IDIP}.  We assign the parameters that fix the mean
$\beta_*$ and the covariance operator $\mathcal{C}_\beta$ as follows:
$\beta_*=1$, $\alpha_\beta=7$, $\bs{\gamma}_\beta=0.01\bs{I}$ and
$\kappa_\beta=0$.  On the far right of Figure \ref{fig: variances}, we
show the resulting spatial variance structure of $\bs{\Gamma}_\beta$
with the weighting. For comparison we also show the typically
implemented case of homogeneous Neumann boundary conditions without
weighting (far left), the case of enforcing Dirichlet boundary
conditions (centre left), and the case of applying a homogeneous Robin
boundary condition following the method of \cite{Roininen2014} (centre
right), We note that the weighted covariance approach nullifies all
boundary effects.
Figure \ref{fig: betadraws} shows three samples drawn from the prior
density, $\mu_\beta$, along with the actual (distributed) Robin
coefficient $\beta_{\rm true}$ used to generate the synthetic data in
both experiments.

\begin{figure}[t!]
\includegraphics[width=\linewidth]{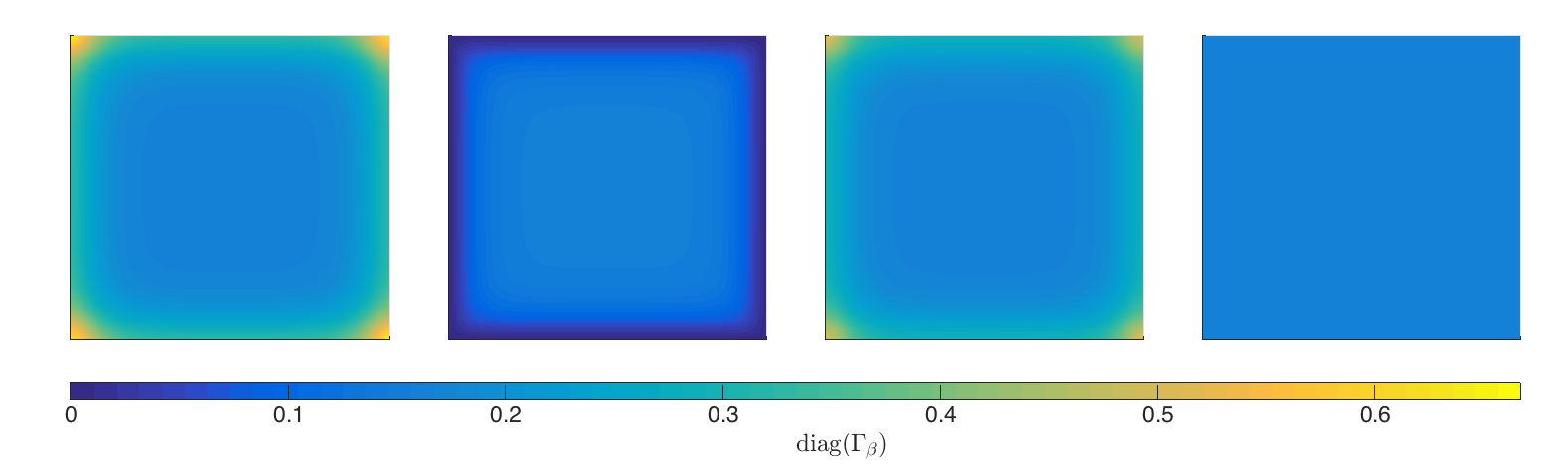}
\caption{The diagonal of the prior covariance operator with
  $\alpha_\beta=7$ and, $\gamma_\beta=0.01$. Far left: With
  homogeneous Neumann boundary conditions. Centre left: With
  homogeneous Dirichlet boundary conditions. Centre right: With
  homogeneous Robin boundary conditions and
  $\kappa_\beta=1.42\sqrt{\gamma_{\beta}/\alpha_{\beta}}$ as in
  \cite{Roininen2014}. Far right: The weighted approach of the current
  paper as discussed in Section \ref{sec: BA2IDIP}, with $\kappa_\beta=0$ (homogeneous Neumann).}\label{fig:
  variances}
\end{figure}

The BAE approach is (in part) based on (approximate) marginalization
over the nuisance parameter.  Technically, this involves drawing
samples from the joint prior density of the nuisance and the primary
parameters $\pi(a,\beta)$ to compute the second order statistics of
$(\bs{\eps},\beta)$ which, in turn, involves the computation of the
forward problem $\bs{f}(a,\beta)$ for the draws.  In this paper, we
take $(a,\beta)$ to be mutually independent.  In addition, we take the
prior on $a(\bs{x})$ to be a Gaussian measure,
$\mu_a=\mathcal{N}(a_*,\mathcal{C}_{a})$ on $L^2(\Omega)$ with
$\mathcal{C}_{a}$ defined similarly as $\mathcal{C}_{\beta}$ in
Section~\ref{sec: BA2IDIP}.  To be precise, we use a squared inverse
elliptic operator as our prior covariance operator with homogeneous
Neumann boundary conditions.  We neglect to weight the covariance
operator $\mathcal{C}_{a}$ by the diagonal of the associated Greens
function and simply impose homogeneous Neumann boundary conditions in
a bid to reduce the computational cost.  Any boundary effects caused
by this are of no consequence as we do not wish to reconstruct $a$.
Formally, the prior covariance operator of $a$ is defined as
$\mathcal{C}_{a}=\mathcal{A}_a^{-2}$, with the operator
$\mathcal{A}_a$ defined (similarly to $\mathcal{C}_\beta$ in Section
\ref{sect: discret}) through the variational problem: For $s\in$
$L^2(\Omega)$, the solution of $\mathcal{A}_aa = s$ satisfies
\begin{align}\label{eq: priorona}
\alpha_a\int_{\Omega}\left(\bs{\gamma}_a\nabla a\cdot\nabla v+a v\right)\;d\bs{x}=\int_{\Omega} s v\;d{\bs x}\quad\text{for all } v\in H^1(\Omega).
\end{align}

\begin{figure}[t!]
\includegraphics[width=\linewidth]{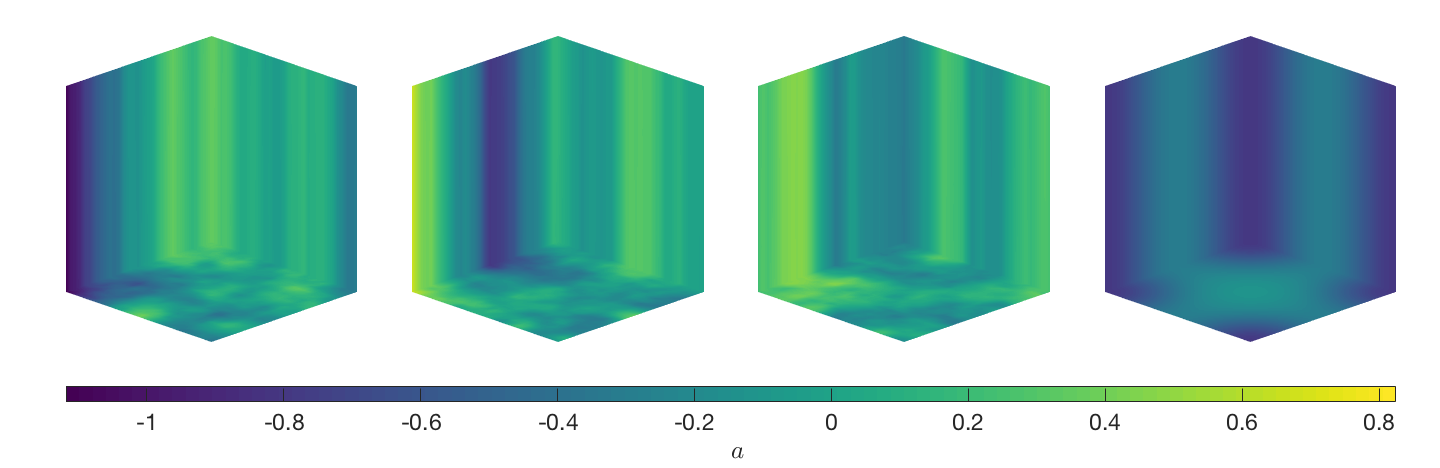}
\caption{ Far left to centre right: Slice plots of samples from the
  prior distribution on $a(\bs{x})$ with isotropic covariance
  structure (first numerical example).  Right: Horizontal cross
  sections of the true value $a_{\rm true}$ used to compute the data.}
\label{fig: sigmadraws}
\end{figure}

\begin{figure}[t!]
\includegraphics[width=\linewidth]{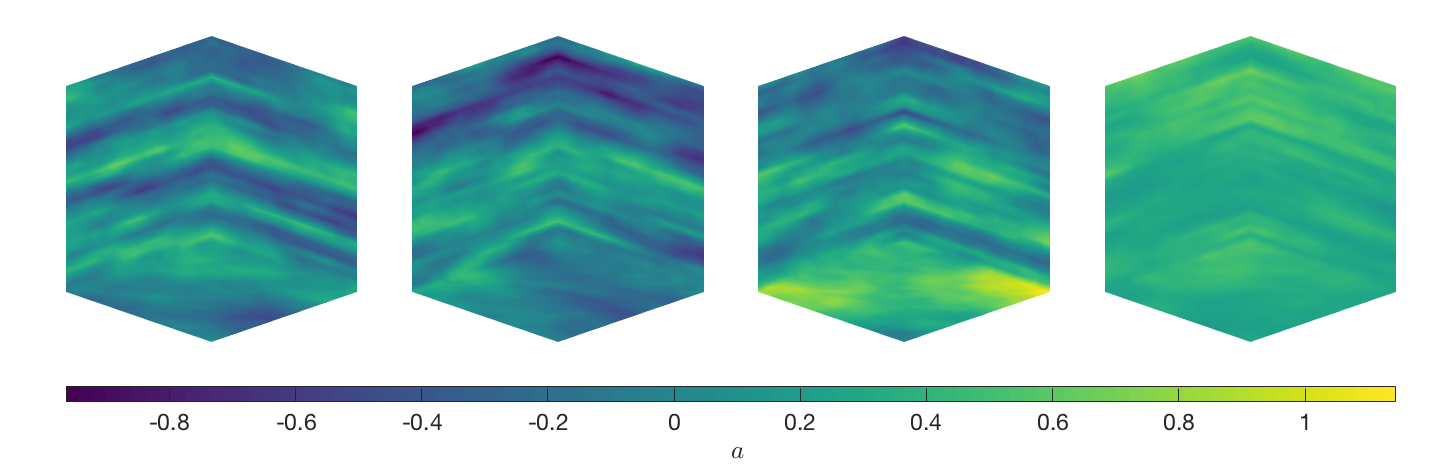}
\caption{ Far left to centre right: Vertical cross sections of samples
  from the prior distribution on $a$ with anisotropic covariance
  structure (second numerical example).  Right: Vertical cross
  sections of the true value $\bs{a}_{\rm true}$ used to compute the
  data.}
\label{fig: Ansigmadraws}
\end{figure}

For the first numerical example with isotropic correlation structure,
we use $\alpha_a=100$ and $\bs{\gamma}_a=10^{-3}\bs{I}$, while for the
second numerical example with anisotropic correlation structure, we
take $\bs{\gamma}_a=\text{diag}(10^{-2},10^{-2},10^{-8})$.  In Figure
\ref{fig: sigmadraws}, three samples drawn from $\mu_a$ for the first
numerical example are shown, along with the true value, $a_{\rm
  true}$, used to generate the synthetic data.  Similarly, in Figure
\ref{fig: Ansigmadraws}, three samples drawn from $\mu_a$ for the
second numerical example are shown along with the true value $a_{\rm
  true}$ used to generate the synthetic data.  A standard requirement
when designing the prior is that the priors should not be ``too
narrow", which here is reflected in the shown draws when compared to
the actual $a(\bs{x})$ and $\beta(\bs{x})$.

\begin{figure}[h!]
\includegraphics[width=\linewidth]{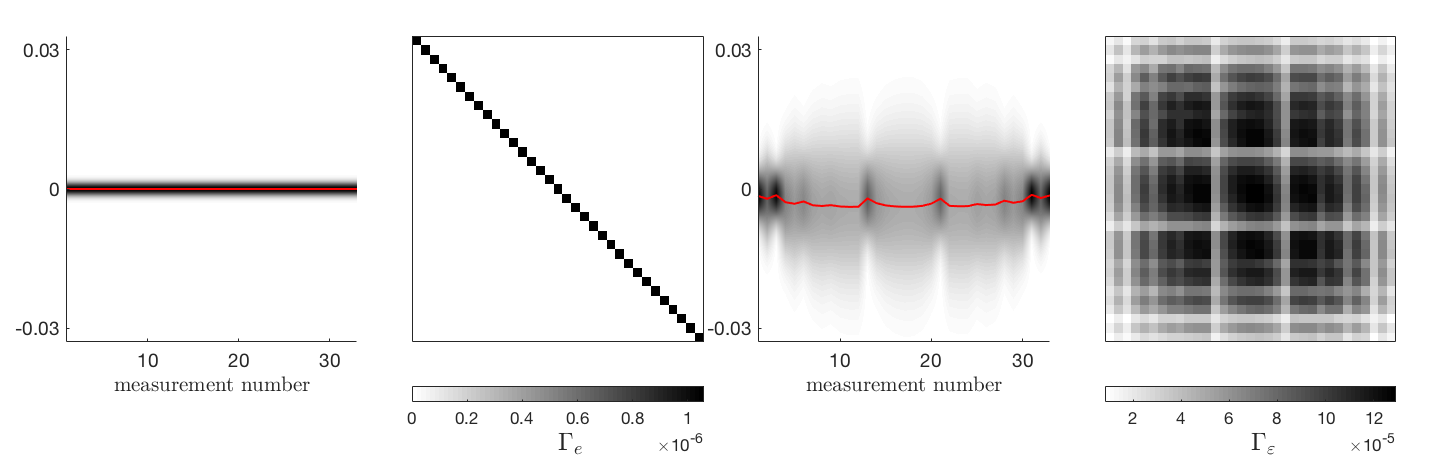}
\caption{Second order statistics of the noise and approximation errors for the second 
(anisotropic conductivity) numerical example. Far left:  The density of the noise, $\mu_{\rm noise}$ (mean shown in red) with higher probability density indicated by darker shading. Centre left: The covariance matrix of the noise $\bs{\Gamma}_{ e}$. Centre right: The density of the total errors, $\mu_{\nu}$ (mean shown in red) with higher probability density indicated by darker shading. Far right: The covariance matrix of the approximation errors $\bs{\Gamma}_{\eps}$.}
\label{fig: Ana_e_stats}
\end{figure}

{\bf Estimation of approximation error statistics.} 
In the linear normal case, that is, ${\bs f}(a,\beta)$ and ${\bs
  f}_{a_*}(\beta)$ both linear and $\pi(a,\beta)$ is normal,
$\bs{\eps}_*$ and $\Gamma_{\bs \eps}$ can be computed analytically.
If this is not the case, both $\bs{\eps}_*$ and $\Gamma_{\bs \eps}$
must be estimated using sample statistics using samples drawn from
(the not necessarily jointly Gaussian) joint prior model
$\pi(a,\beta)$. With an ensemble of $r$ samples,
$(\beta^{(\ell)},a^{(\ell)})$ from the associated prior densities, we
compute
\begin{align}
{\bs \eps}^{(\ell)}={\bs f}(a^{(\ell)},\beta^{(\ell)})-{\bs f}_{a_*}(\beta^{(\ell)}),\quad \ell=1,2,,\dots,r,
\end{align}
and take the mean and covariance as
\begin{align}
\bs{\eps}_*=\frac{1}{r}\sum_{\ell=1}^r \bs{\eps}^{(\ell)}\quad\text{and}\quad\bs{\Gamma}_{\eps}=\frac{1}{1-r}\sum_{\ell=1}^r (\bs{\eps}^{(\ell)}-\bs{\eps}_*)(\bs{\eps}^{(\ell)}-\bs{\eps}_*)^T.
\end{align}
The number of samples $r$ required depends on the models, the variance
of the approximation error, and the joint prior model, see for example
\cite{Kaipio2007a}. However, we remark that all samples (and all
accurate forward simulations) are carried out at the offline stage and
that the accurate forward model is never used in the inversion, with
only the approximate model evaluated at the online stage.  For the
current problem, 1000 samples were drawn for both numerical examples
to compute the approximation error statistics. Figures \ref{fig:
  a_e_stats} and \ref{fig: Ana_e_stats} show the statistics of the
noise and of the approximation errors for the first and second
examples, respectively.  It is clear that the approximation errors
dominate the noise in both cases, entries on the diagonal of
$\bs{\Gamma}_{\eps}$ being almost two orders of magnitude larger than
those on the diagonal of $\bs{\Gamma}_e$. Thus, not including the
approximation errors will most likely lead to meaningless results, as
confirmed below in Section~\ref{sec:Results}.  Note the difference in
the mean and covariance of $\bs{\eps}$ due to different priors on $a$.

\begin{figure}[t!]
\includegraphics[width=\linewidth]{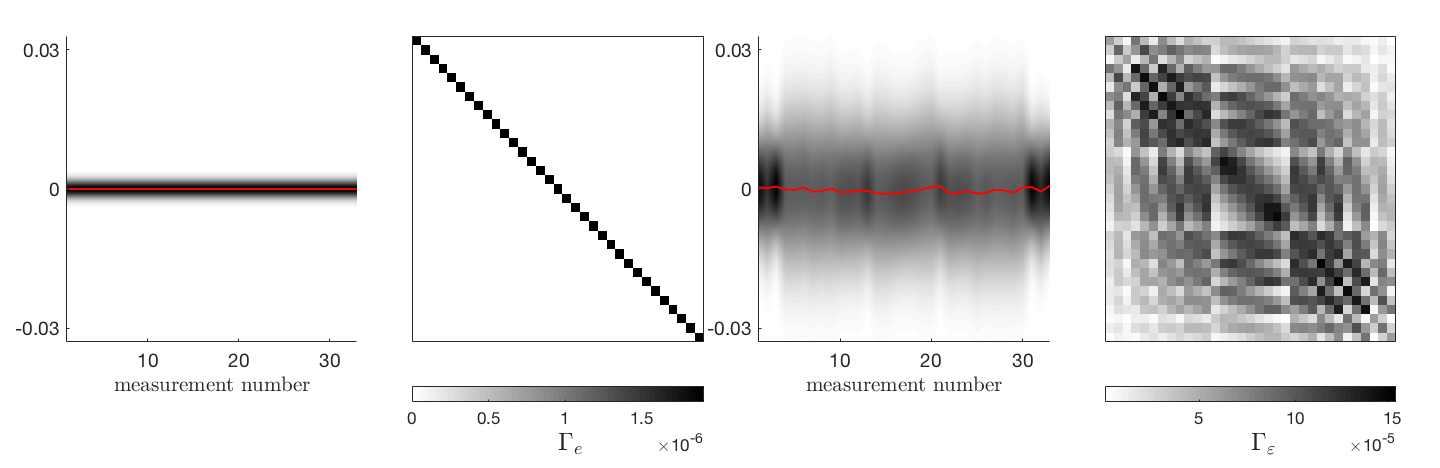}
\caption{Second order statistics of the noise and approximation errors for the first (isotropic conductivity) numerical example. Far left:  The density of the noise, $\mu_{\rm noise}$ (mean shown in red) with higher probability density indicated by darker shading. Centre left: The covariance matrix of the noise $\bs{\Gamma}_{\rm e}$. Centre right: The density of the total errors, $\mu_{\nu}$ (mean shown in red). Far right: The covariance matrix of the approximation errors $\bs{\Gamma}_{\eps}$.}
\label{fig: a_e_stats}
\end{figure}

{\bf The MAP and approximate posterior covariance estimates.} 
To compare the solutions of the inverse problems formulated with the
approximation error noise model and with the conventional error model, we
compute the respective MAP estimates and posterior covariances
matrices. For both the isotropic and anisotropic conductivity cases,
we compute the following three MAP estimates:
\vspace{-0.1in}
\begin{itemize}
\item MAP-REF: The reference maximum a posteriori estimate for $\beta$ with the conventional error model and use of the {\em actual} value of $a$ in the model, that is, $ \bs {d}^{\rm obs}=\bs{f}({a_{\rm true}},\beta)+{\bs e}=\bs{f}_{a_{\rm true}}(\beta)+{\bs e}$. This reconstruction is taken as the benchmark one, as it is computed with no modelling errors present. The estimate is computed as
\begin{align}\label{eq: mapref_est}
\bs{\beta}^{\rm REF}_{\rm MAP}=\min_{\bs{\beta}\in \mathbb{R}^{n}}\frac{1}{2}\normm[\bs{\Gamma}_{\bs e}^{-1}]{\bs{f}_{\bs{a}_{\rm true}}(\bs{\beta})-\bs{d}^{\rm obs}} ^2 +\frac{1}{2}\normm[\bold{M}]{\bold{A}\left(\boldsymbol{\beta}-\boldsymbol{\beta}_*\right)}^2.
\end{align}

\item MAP-CEM: The  maximum a posteriori estimate for $\beta$ with the conventional error model, using the (incorrect) fixed $a=a_*$ in the forward model, that is, $ \bs {d}^{\rm obs}=\bs{f}_{a_*}(\beta)+{\bs e}$. This estimate is computed as
\begin{align}\label{eq: mapcem_est}
\bs{\beta}^{\rm CEM}_{\rm MAP}=\min_{\bs{\beta}\in \mathbb{R}^{n}}\frac{1}{2}\normm[\bs{\Gamma}_{\bs e}^{-1}]{\bs{f}_{a_*}(\bs{\beta})-\bs{d}^{\rm obs}} ^2 +\frac{1}{2}\normm[\bold{M}]{\bold{A}_\beta\left(\boldsymbol{\beta}-\boldsymbol{\beta}_*\right)}^2.
\end{align}

\item MAP-BAE: The  maximum a posteriori estimate for $\beta$ with the approximation error model and using fixed $a=a_*$ in the model, that is, $ \bs {d}^{\rm obs}=\bs{f}_{a_*}(\beta)+{\bs \nu}$. This estimate is computed as
\begin{align}\label{eq: mapaem_est}
\bs{\beta}^{\rm BAE}_{\rm MAP}=\min_{\bs{\beta}\in \mathbb{R}^{n}}\frac{1}{2}\normm[\bs{\Gamma}_{\bs \nu}^{-1}]{\bs{f}_{a_*}(\bs{\beta})-\bs{d}^{\rm obs}+\bs{\nu}_* } ^2 +\frac{1}{2}\normm[\bold{M}]{\bold{A}_\beta\left(\boldsymbol{\beta}-\boldsymbol{\beta}_*\right)}^2.
\end{align}

\end{itemize}

The related approximate posterior covariance matrices are then as follows.

\begin{itemize}
\item $\Gamma^{\rm REF}_{\rm post}$: The reference posterior
  covariance matrix is computed using the conventional error model and
  using the actual value of $a$ in the model, i.e. $ \bs {d}^{\rm
    obs}=\bs{f}(a_{\rm true},\beta)+{\bs e}$. Since, in this case,
  there are no modelling errors present, we expect the reference
  posterior covariance matrix to be smaller (in the sense of quadratic
  forms) than the posterior covariance matrix obtained using the
  approximation error model. This posterior covariance matrix is
\begin{align}\label{eq: covref_est}
\Gamma^{\rm REF}_{\rm post}=\left(\bs{F}_{a_{\rm true}}^\natural\bs{\Gamma}_{\bs e}^{-1}\bs{F}_{a_{\rm true}}+\bs{\Gamma}_{\bs{\beta}}^{-1}\right)^{-1}.
\end{align}

\item $\Gamma^{\rm CEM}_{\rm post}$: The posterior covariance matrix
  with the conventional error model is computed using the conventional
  error model with the fixed $a=a_*$ in the model, that is, $ \bs
  {d}^{\rm obs}=\bs{f}_{a_*}(\beta)+{\bs e}$. This posterior
  covariance matrix is
\begin{align}\label{eq: covcem_est}
\Gamma^{\rm CEM}_{\rm post}=\left(\bs{F}_{a_*}^\natural\bs{\Gamma}_{\bs e}^{-1}\bs{F}_{a_*}+\bs{\Gamma}_{\bs{\beta}}^{-1}\right)^{-1}.
\end{align}

\item $\Gamma^{\rm BAE}_{\rm post}$: The approximation error model
  posterior covariance matrix is computed using the approximation
  error model and using fixed $a=a_*$ in the model, i.e., $ \bs
  {d}^{\rm obs}=\bs{f}_{a_*}(\beta)+{\bs \nu}$. This posterior
  covariance matrix is
\begin{align}\label{eq: covaem_est}
\Gamma^{\rm BAE}_{\rm post}=\left(\bs{F}_{a_*}^\natural\bs{\Gamma}_{\bs \nu}^{-1}\bs{F}_{a_*}+\bs{\Gamma}_{\bs{\beta}}^{-1}\right)^{-1}.
\end{align}

\end{itemize}

\subsection{Results}
\label{sec:Results}
\dontshow{The primary goal of this section is to compare the
  performance of the conventional error model approach with the
  Bayesian approximation error approach when we take the simplifying
  approximation that the conductivity is approximated by its mean,
  i.e., $a=a_*$, for the inverse problem outlined in section
  \ref{sect: RecRC}. For comparison we also include the reference
  results, i.e., the results obtained when we use the correct value
  for $a$, $a=a_{\rm true}$ in the reconstructions. Initially, we
  discuss the MAP point estimates outlined in equations~\eqref{eq:
    mapref_est}-\eqref{eq: mapaem_est} and then analyze
  the posterior covariance matrices, and investigate the
  approximation of the prior-preconditioned data misfit Hessians
  outlined in section~\ref{section: MapAppPost}.}

The computation of the reference MAP estimate, $\bs{\beta}^{\rm
  REF}_{\rm MAP}$ given in (\ref{eq: mapref_est}), the conventional
error model MAP estimate, $\bs{\beta}^{\rm CEM}_{\rm MAP}$ given in
(\ref{eq: mapcem_est}), and the approximation error model MAP
estimate, $\bs{\beta}^{\rm BAE}_{\rm MAP}$ given in (\ref{eq:
  mapaem_est}), is done by applying an inexact adjoint-based
Gauss-Newton method outlined in Section \ref{section: adjointopt}. We
start each of the optimization procedures with the prior mean as the
initial guess, that is $\bs{\beta}_0=\bs{\beta}_*=1$. A preconditioned
conjugate gradient (CG) method is used with an Einsentat-Walker
condition \cite{Eisenstat1996} which terminates the CG iterations
early, when the norm of the gradient is sufficiently reduced. In line
with \cite{Bui-Thanh2013}, we use the prior operator as a
preconditioner for the CG iterations.

\begin{figure}[h!]
\includegraphics[width=\linewidth]{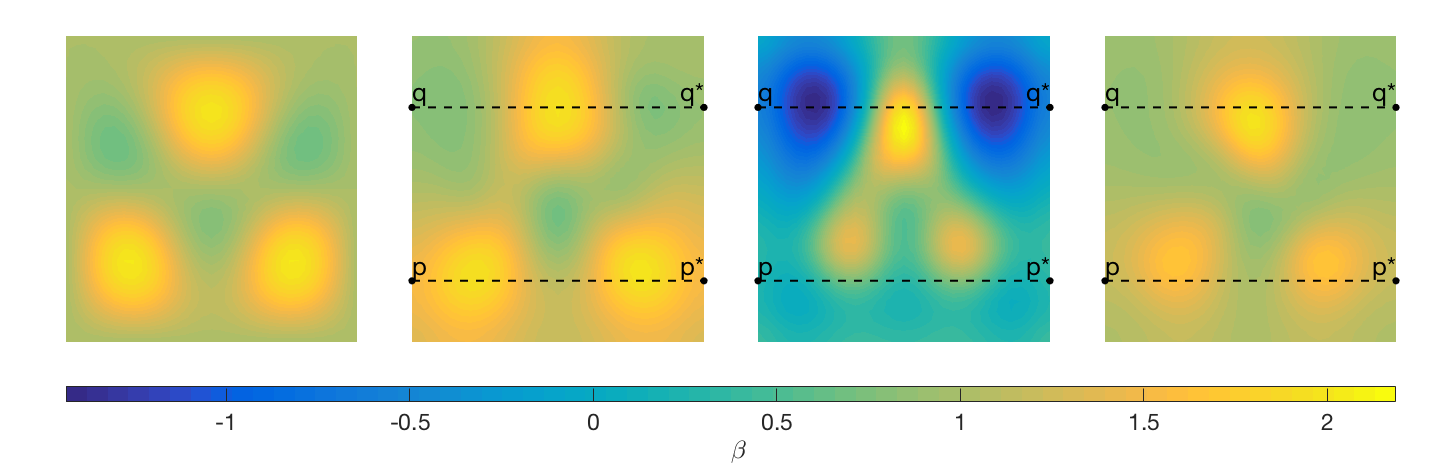}
\caption{Isotropic conductivity. Far left: The true Robin coefficient $\beta_{\rm
    true}$. Centre left to far right: The reference MAP estimate
  $\beta_{\rm MAP}^{\rm REF}$, the conventional error model MAP
  estimate $\beta_{\rm MAP}^{\rm CEM}$ and the approximation error
  model MAP estimate $\beta_{\rm MAP}^{\rm BAE}$, respectively.  The
  cross sections are shown in
  Figures~\ref{fig:pt1_recons}-\ref{fig:pt2_recons}.}
\label{fig:estimates}
\end{figure}

\begin{figure}[h!]
\includegraphics[width=\linewidth]{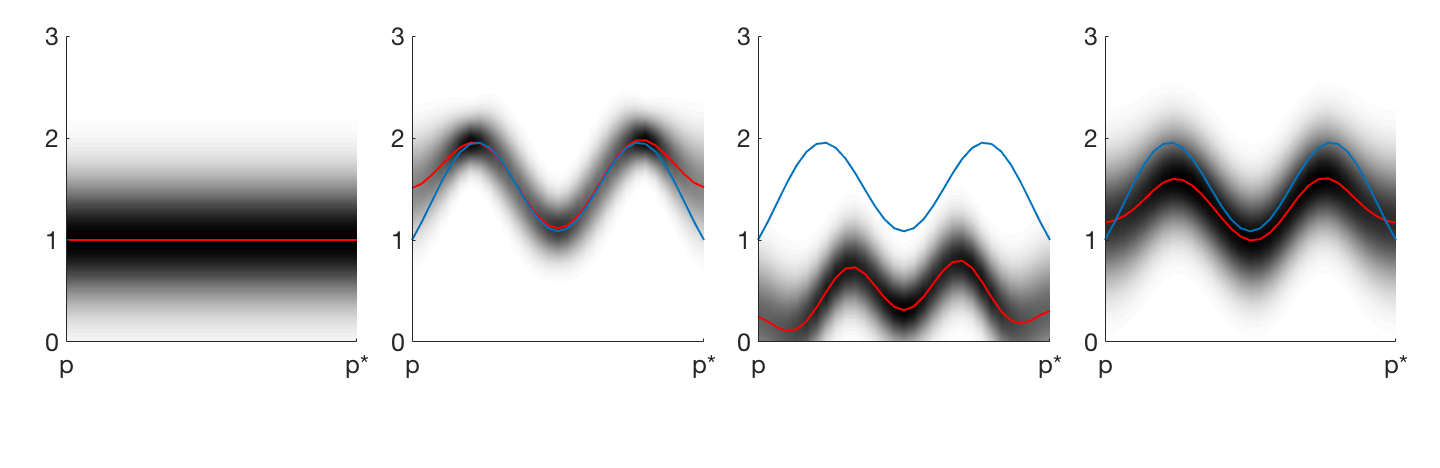}
\caption{Isotropic conductivity.  Far left: the prior mean and
  marginal distibutions of the $\beta$ along the line p-p$^\ast$.
  Centre left to far right: The cross sections of Figure
  \ref{fig:estimates} along the line p-p$^\ast$.  The reference MAP
  estimate $\beta_{\rm MAP}^{\rm REF}$, the conventional error model
  MAP estimate $\beta_{\rm MAP}^{\rm CEM}$ and the approximation error
  model MAP estimate $\beta_{\rm MAP}^{\rm BAE}$.  The true $\beta$
  and the MAP estimates are shown in blue and red, respectively, along
  with the approximate posterior marginal distributions of
  $\beta$. }
\label{fig:pt1_recons}
\end{figure}

\begin{figure}[h!]
\includegraphics[width=\linewidth]{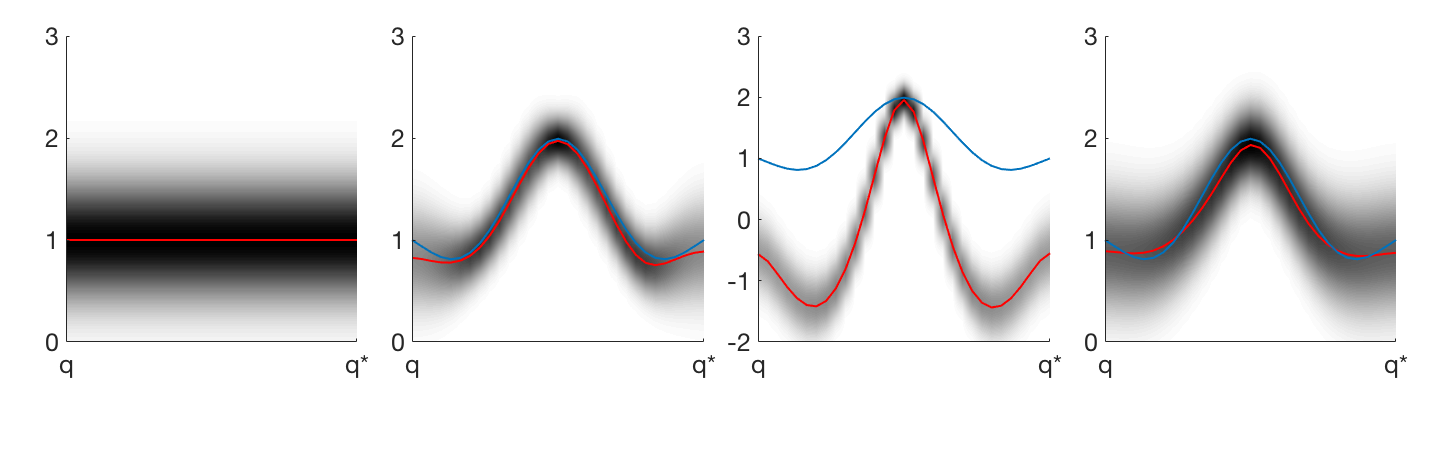}
\caption{Isotropic conductivity.
The cross sections of Figure \ref{fig:estimates} along the line q-q$^\ast$.
The discription otherwise as in Figure~\ref{fig:pt1_recons}.}
\label{fig:pt2_recons}
\end{figure}

{\bf The estimates with isotropic conductivity.} 
On the far left of Figure \ref{fig:estimates} we show the true basal
Robin coefficient, $\beta_{\rm true}$, which is used to generate the
measurements. We also show the reconstructed reference MAP estimate
(centre left), the reconstructed conventional error model MAP estimate
(centre right), and the reconstructed approximation error model MAP
estimate (far right). The images of the reconstructions in Figure
\ref{fig:estimates} also show the (dotted) lines p-p$^*$ and q-q$^*$
which are the locations of the cross sections shown in Figures
\ref{fig:pt1_recons} and \ref{fig:pt2_recons}.  We now discuss several 
observations that can be made  from Figures  \ref{fig:pt1_recons} and
 \ref{fig:pt2_recons}. Firstly, the reference estimate
is clearly {\it feasible} in the sense that the posterior uncertainty
supports the actual Robin coefficient. Conversely, the estimates with the
conventional error model (severely) underestimate the true basal Robin
coefficient.  In particular, the estimate is clearly {\it
  infeasible}: The actual coefficient has almost vanishing posterior
density at almost all points along the cross sections. On the other
hand, the estimate with the approximation error model is clearly a
feasible one with the posterior marginals supporting the actual Robin
coefficient.  Finally, we see that marginalization over the conductivity
results in the widening of the posterior density which is evident when
comparing the marginal densities of the reference and approximation
error estimates.

\begin{figure}[h!]
\includegraphics[width=\linewidth]{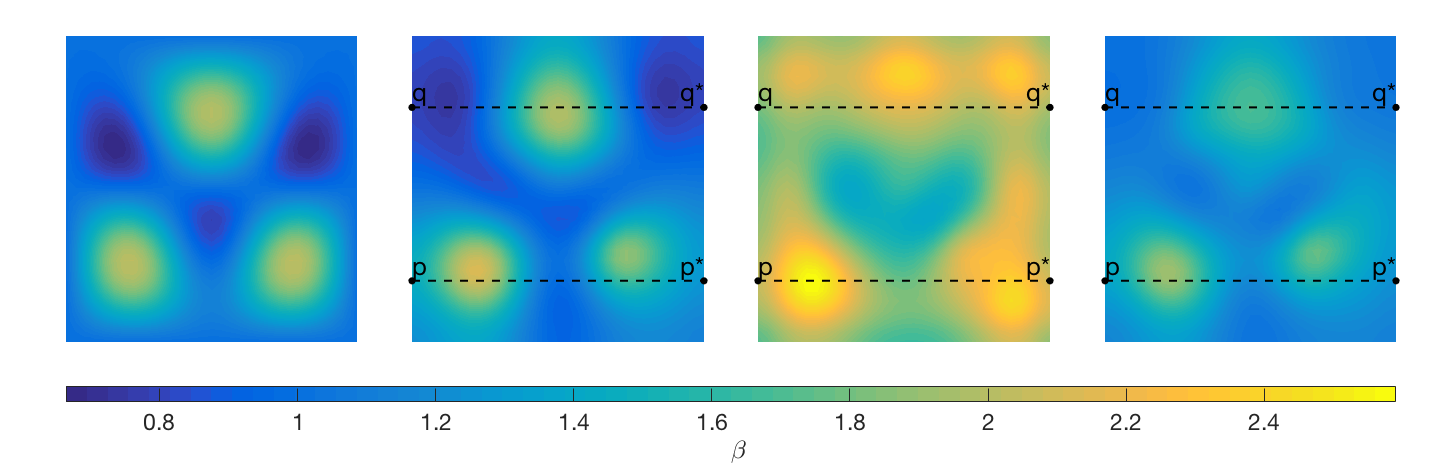}
\caption{Anisotropic conductivity.  Far left: The true Robin
  coefficient $\beta_{\rm true}$. Centre left to far right: The
  reference MAP estimate $\beta_{\rm MAP}^{\rm REF}$, the conventional
  error model MAP estimate $\beta_{\rm MAP}^{\rm CEM}$ and the
  approximation error model MAP estimate $\beta_{\rm MAP}^{\rm BAE}$,
  respectively.  The cross sections are shown in
  Figs.~\ref{fig:Anpt1_recons}-\ref{fig:Anpt2_recons}.}
\label{fig:Anestimates}
\end{figure}

\begin{figure}[h!]
\includegraphics[width=\linewidth]{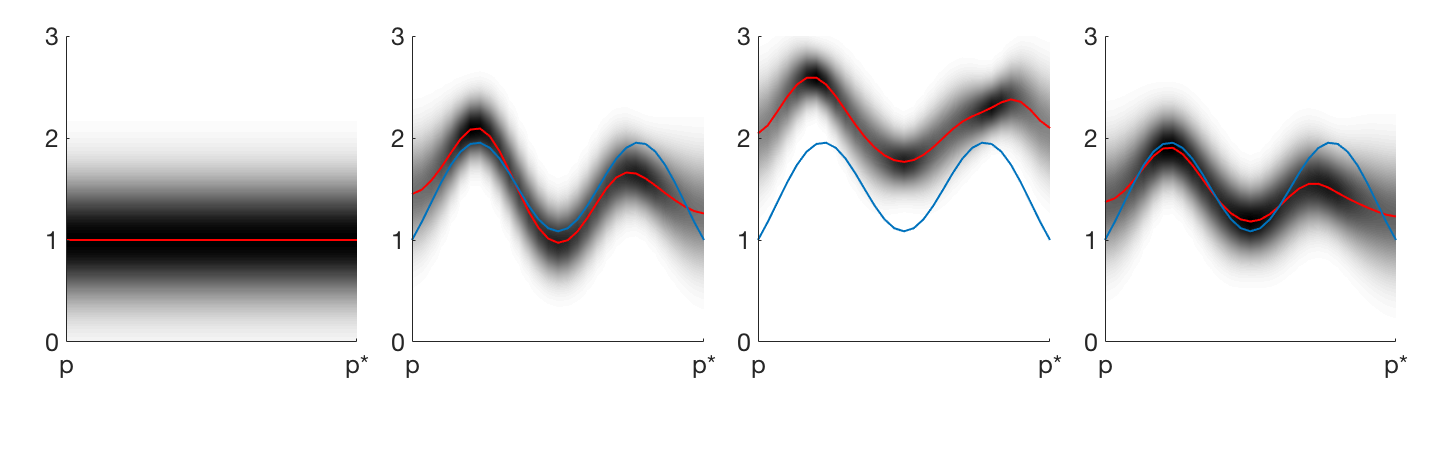}
\caption{Anisotropic conductivity.  The cross sections of Figure
  \ref{fig:Anestimates} along the line q-q$^\ast$.  The discription
  otherwise as in Figure~\ref{fig:pt1_recons}. }
\label{fig:Anpt1_recons}
\end{figure}

\begin{figure}[h!]
\includegraphics[width=\linewidth]{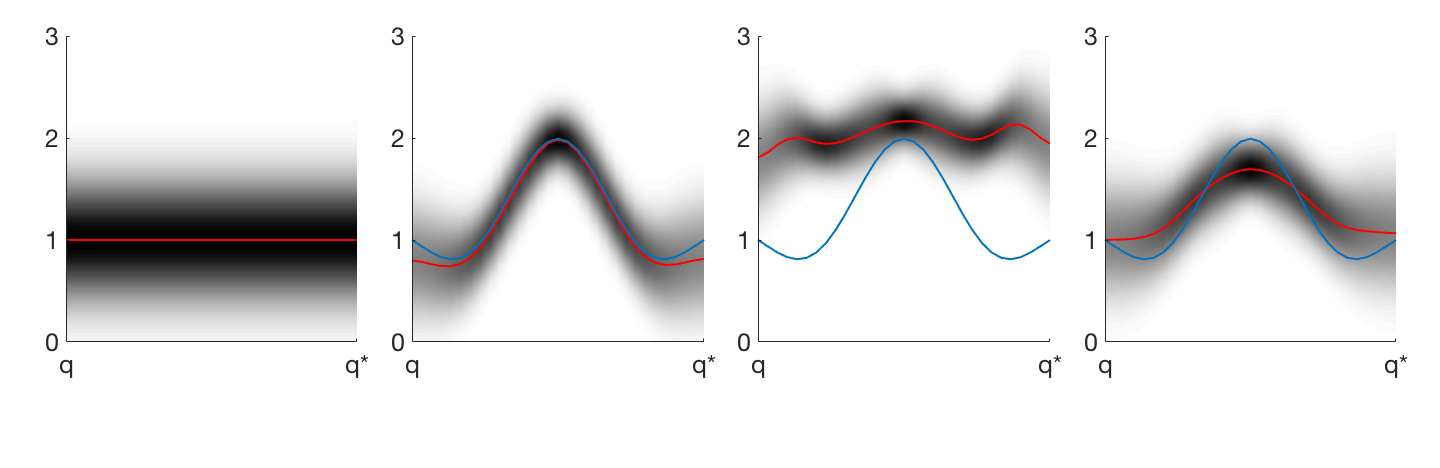}
\caption{Anisotropic conductivity.
The cross sections of Figure \ref{fig:Anestimates} along the line q-q$^\ast$.
The discription otherwise as in Figure \ref{fig:pt1_recons}.}
\label{fig:Anpt2_recons}
\end{figure}

{\bf The estimates with anisotropic conductivity.}  The corresponding
results in the case of the anisotropic conductivity are shown in
Figures~\ref{fig:Anestimates}-\ref{fig:Anpt2_recons}.  The results are
qualitatively similar to the isotropic case.  In this case, the
conventional error severely overestimates the actual Robin
coefficient, the only difference between the two cases being the
spatial covariance structure of the conductivity $a$.  With the
approximation error model, the estimates are still feasible but
slightly worse when comparing to the isotropic conductivity case.

{\bf Computational costs.}
In this section, we compare the computational cost of the inverse
solution method applied for the three methods: the reference case; the conventional error procedure; and the BAE approach. We measure this cost
in terms of number of Poisson solves needed for the optimization
algorithm to converge. We note that to compute the MAP estimates, the
number of Poisson solves needed per Gauss-Newton iteration can be
calculated as
{\fontseries{b}\selectfont\#}{Poisson}=2+2{\fontseries{b}\selectfont\#}{CG}+{\fontseries{b}\selectfont\#}{back.}
where {\fontseries{b}\selectfont\#}{CG} is the number of CG iterations
and {\fontseries{b}\selectfont\#}{back.} is the number of back-tracks
needed to get a sufficient decrease in the objective function. In all
three cases, the convergence of the Gauss-Newton iterations is
established when the norm of the gradient (relative to the initial
norm of the gradient) is decreased by a factor of $10^7$. The results
shown in Table~\ref{table: costs} reveal that, at the inversion stage,
the BAE approach is approximately as expensive as the conventional
error approach. For completeness, we also show the costs in the
reference reconstructions.

\begin{table}
\caption{Comparison of the computational costs.  The cost of solving
  for the MAP estimates in the reference case (REF), the case of using
  the conventional error model  (CEM) with fixed $a=a_*$ and the
  case of using the approximation error model (BAE) with fixed $a=a_*$, measured in
  number of Poisson solves. The first column ({\bf MAP}) refers to
  which MAP estimate we are solving for, MAP-REF (\ref{eq:
    mapref_est}), MAP-CEM (\ref{eq: mapcem_est}), or MAP-BAE (\ref{eq:
    mapaem_est}); the second column ({\fontseries{b}\selectfont\#}{\bf
    GN}) reports the number of Gauss-Newton iterations; the third
  ({\fontseries{b}\selectfont\#}{\bf CG}) and fourth ({\bf avg.CG})
  columns show total and the average (per Gauss-Newton iteration)
  number of CG iterations; the fifth column
  ({\fontseries{b}\selectfont\#}{\bf back}) reports the total number
  of backtracks needed throughout the Gauss-Newton iterations; and the
  last column ({\fontseries{b}\selectfont\#}{\bf Poisson}) reports the
  total number of Poisson solves (from forward, adjoint, and
  incremental forward and adjoint problems). The Gauss-Newton
  iterations are terminated when the norm of the gradient is decreased
  by a factor of $10^7$, while the CG iterations are terminated inline
  with the Einstat-Walker condition \cite{Eisenstat1996}. These
  results illustrate that the use of the approximation error approach can
  be carried out at no additional cost compared to the conventional
  error approach and reference case.}\label{table: costs}
\noindent
\begin{tabu} to \textwidth {lXXXXXr}
   \toprule
   & {\bf MAP} & {\fontseries{b}\selectfont\#}{\bf GN} &{\fontseries{b}\selectfont\#}{\bf CG} & {\bf avg.CG} &  {\fontseries{b}\selectfont\#}{\bf back} & {\fontseries{b}\selectfont\#}{\bf Poisson}\\
   \bottomrule
    Example 1 & & & & & &\\
  & \myalign{l}{REF} & \myalign{l}{8} & \myalign{l}{117} & \myalign{l}{15} & \myalign{l}{0} & \myalign{l}{250}\\
  & \myalign{l}{CEM} & \myalign{l}{11} & \myalign{l}{101} & \myalign{l}{10} & \myalign{l}{4} & \myalign{l}{228}\\
  & \myalign{l}{BAE} & \myalign{l}{5} & \myalign{l}{57} & \myalign{l}{12} & \myalign{l}{0} & \myalign{l}{124}\\
     \bottomrule
    Example 2 & & & & & &\\
  & \myalign{l}{REF} & \myalign{l}{6} & \myalign{l}{54} & \myalign{l}{9} & \myalign{l}{0} & \myalign{l}{120}\\
  & \myalign{l}{CEM} & \myalign{l}{8} & \myalign{l}{95} & \myalign{l}{12} & \myalign{l}{0} & \myalign{l}{206}\\
  & \myalign{l}{BAE} & \myalign{l}{5} & \myalign{l}{97} & \myalign{l}{20} & \myalign{l}{0} & \myalign{l}{204}\\
   \bottomrule
\end{tabu}
\end{table}

{\bf Interpretation of the posterior covariance matrices.}  Here we
discuss and compare the three posterior covariance matrices
corresponding to the reference case, the conventional error with
reduced model case, and the approximation error model case, which are
defined in (\ref{eq: covref_est})--(\ref{eq: covaem_est}). We begin
the discussion by analyzing the spectrum of the respective
prior-preconditioned data misfit Hessian components of the posterior
covariance matrices for both numerical examples.
 
\begin{figure}[t!]
\centering
\includegraphics[width=0.4\linewidth]{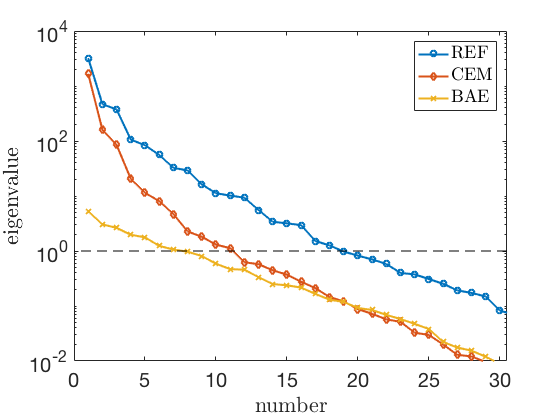}  \includegraphics[width=0.4\linewidth]{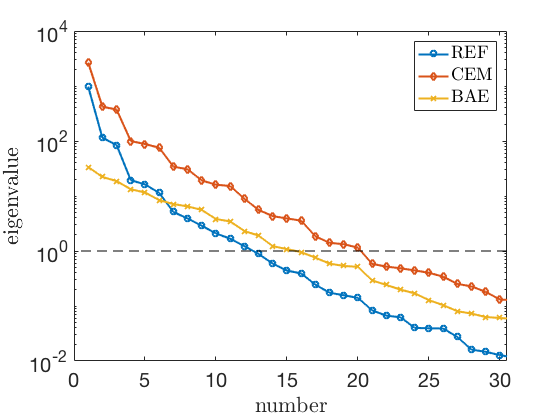}
\caption{Semi-log plots of the eigenvalues of the prior-preconditioned
  misfit Hessians for the isotropic case (left) and the anisotropic
  case (right).}

\label{fig: ppHmisspect}
\end{figure}

On the left of Figure \ref{fig: ppHmisspect} is shown the dominant
spectrum of the prior-preconditioned data misfit Hessian for the three
cases evaluated at the respective MAP estimates for the first example,
while the same results for the second example are shown on the right
(reference case with blue circles, conventional error with reduced
model with red diamonds, and approximation error model with yellow
crosses). In all three cases we are only required to retain relatively
few eigenvalues to compute a reasonable low rank approximation of the
Hessian. Specifically, in the first example, for the reference case we
need about 30 eigenvalues, while for both the conventional error model
with reduced model and the approximation error model we need bout
20. In the second example the reference case requires the retention of
about 20 eigenvalues, the conventional error model with reduced model
requires about 30, and the reduced model with approximation error
model we need 25. We note that that the these low numerical ranks are
substantially smaller than the 961 degrees of freedom of the parameter
(i.e., we see a compression of the parameter dimension of about
30). Hence the approximate posterior covariance matrix along with draws
from the posterior can be cheaply computed using (\ref{eq: svdinvH}) and
 (\ref{eq: SQRTsvdinvH})
 respectively.

\begin{figure}[h!]
\includegraphics[width=\linewidth]{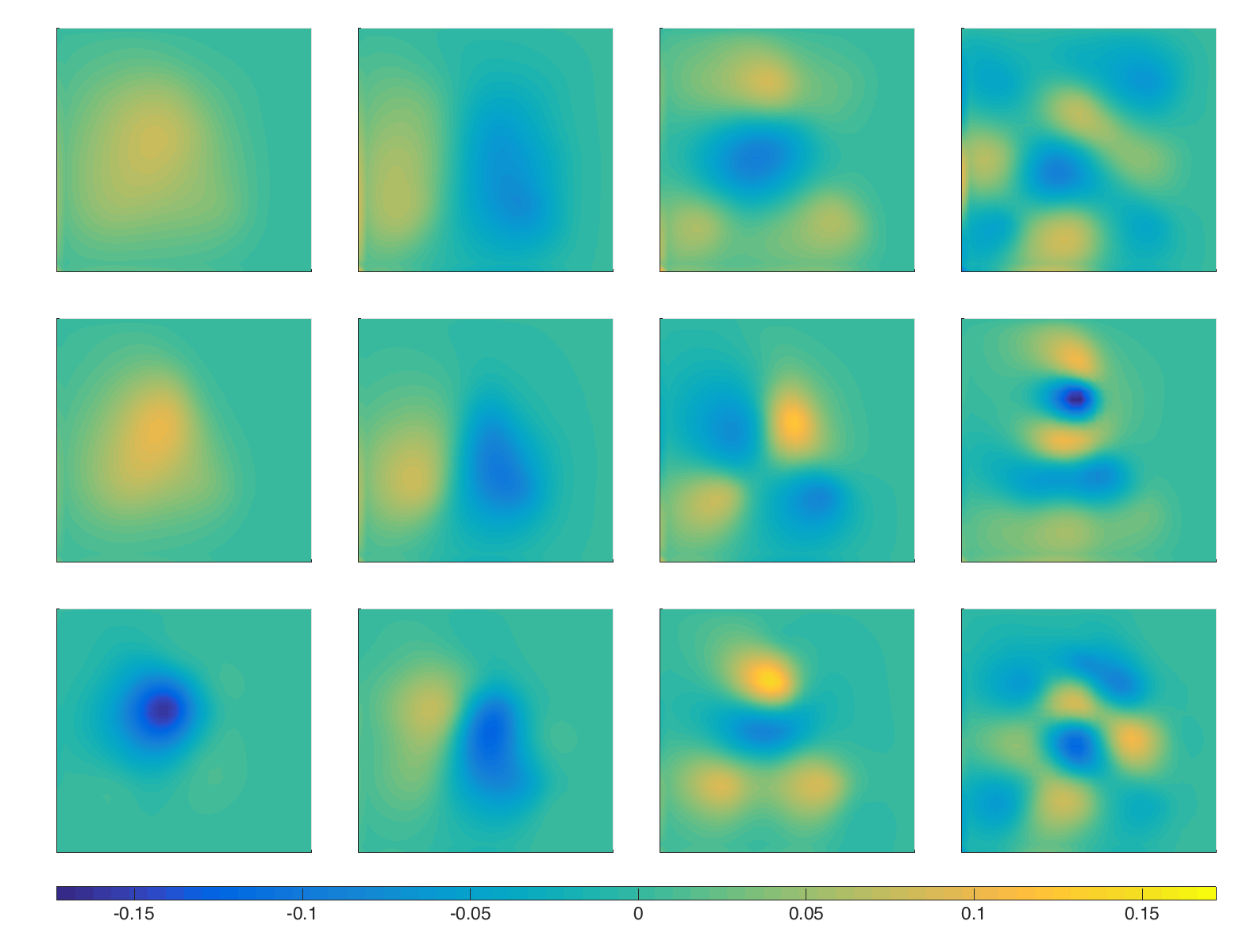}
\caption{Eigenvectors of the reference prior-preconditioned misfit
  Hessian (top row), conventional error model prior-preconditioned
  misfit Hessian (middle row), and approximation error model
  prior-preconditioned misfit Hessian (bottom row) computed at
  $\bs{\beta}_{\rm MAP}^{\rm REF}$, $\bs{\beta}_{\rm MAP}^{\rm CEM}$,
  and $\bs{\beta}_{\rm MAP}^{\rm BAE}$ respectively. From left to
  right: The eigenvectors corresponding to the first (i.e., the
  largest), the third, the fifth and the tenth eigenvalues.}
\label{fig: ppHmvects}
\end{figure}

\begin{figure}[h!]
\includegraphics[width=\linewidth]{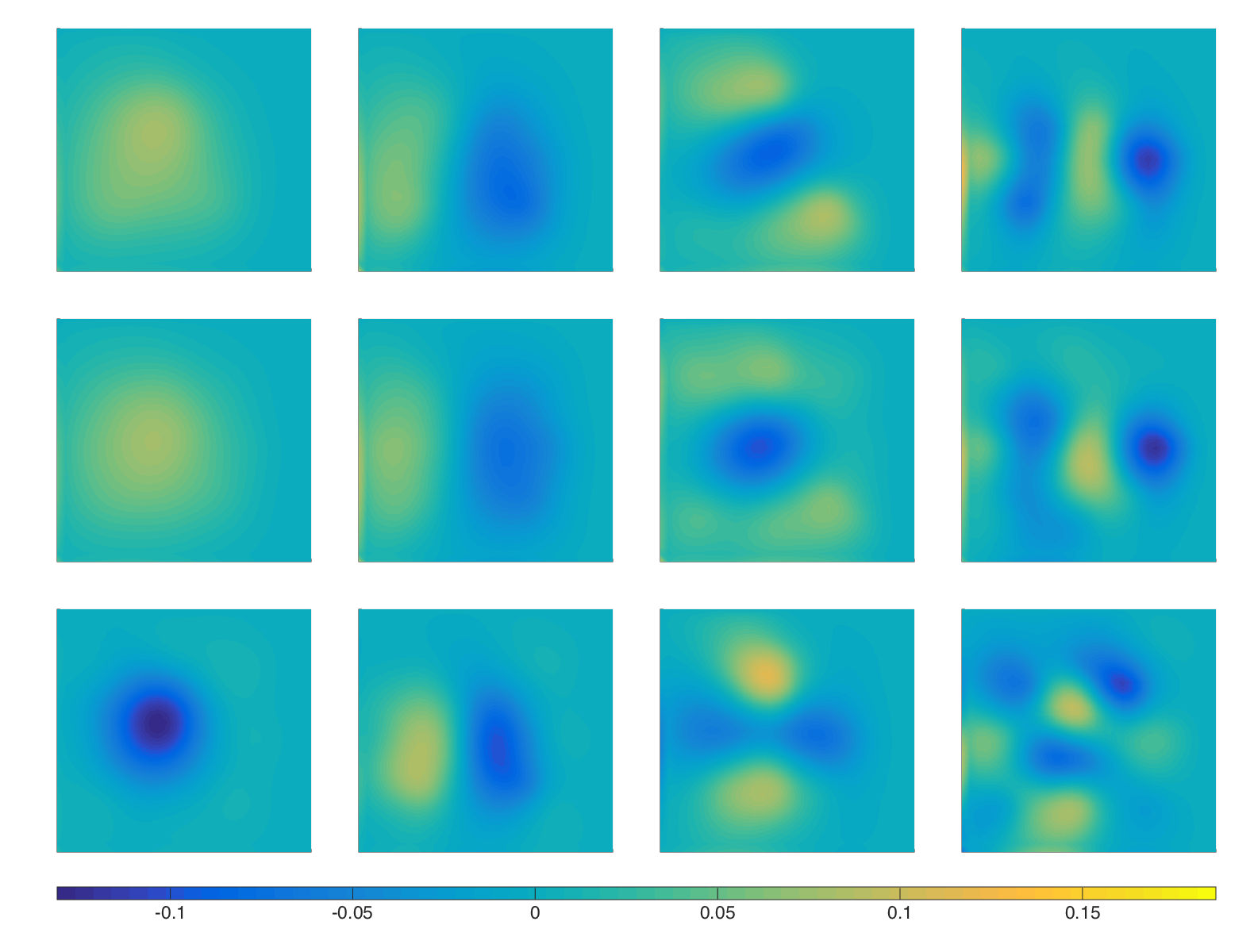}
\caption{Eigenvectors of the
  reference prior-preconditioned misfit Hessian (top row),
  conventional error model prior-preconditioned misfit Hessian (middle
  row), and approximation error model prior-preconditioned misfit
  Hessian (bottom row) computed at $\bs{\beta}_{\rm MAP}^{\rm REF}$,
  $\bs{\beta}_{\rm MAP}^{\rm CEM}$, and $\bs{\beta}_{\rm MAP}^{\rm
    BAE}$ respectively. From left to right: The eigenvectors
  corresponding to the first (i.e., the largest), the third, the fifth
  and the tenth eigenvalues.}
\label{fig: AnppHmvects}
\end{figure}

In Figures \ref{fig: ppHmvects} and \ref{fig: AnppHmvects} we show
several eigenvectors corresponding to the dominant eigenvalues of the
prior-preconditioned data misfit Hessians corresponding to the three
cases. The most dominant eigenvalues can be interpreted as the modes
in the basal Robin coefficient for which the data contains the most
information about. The first few eigenvectors of all three cases are
fairly similar (note the sign change for the first eigenvector is the
approximation error model case). However as the level of oscillation
in the eigenqvectors increases we see that the differences between
corresponding eigenvectors between the three models increase.


\begin{figure}[h!]
\includegraphics[width=\linewidth]{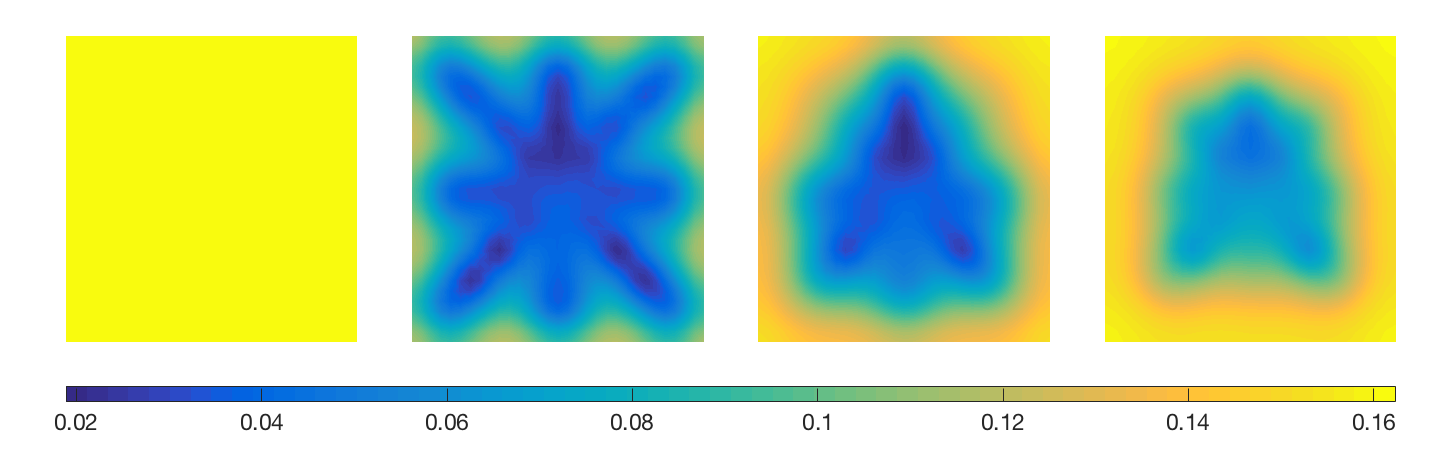}
\caption{The reduction in variance for the isotropic case. Far left:
  The diagonal of the prior covariance matrix, $\bs{\Gamma}_\beta$, as
  outlined in section \ref{sec: BA2IDIP}. Centre left: The diagonal of
  the reference approximate posterior covariance matrix,
  $\bs{\Gamma}^{\rm REF}_{\rm post}$. Centre right: The diagonal of
  the conventional error model approximate posterior covariance
  matrix, $\bs{\Gamma}^{\rm CEM}_{\rm post}$. Far right: The diagonal
  of the approximation error model approximate posterior covariance
  matrix, $\bs{\Gamma}^{\rm BAE}_{\rm post}$.}\label{fig: var_reduc}
\end{figure}

We can further asses the uncertainty in the estimates by analyzing the
full pointwise posterior variances, i.e. the diagonal of the posterior
covariance matrices. Such analysis also provides insight into how much
the variance (from the prior to the posterior) is reduced by taking
the data into consideration. In Figures \ref{fig: var_reduc} and
\ref{fig: var_reduc_an} we show the diagonal of the prior covariance
matrix and the three posterior covariance matrices, for the first and
second numerical examples, respectively.

Most information, and consequently the greatest reduction in
uncertainty, is directly below the measurement locations and in areas
where the parameter field attains relatively high values. The
reduction in variance in the vicinity of the measurement locations is
typical for inverse problems. On the other hand, it is evident that
the reduction in variance is greatest where the parameter achieves
higher values. This is due to the fact that for higher values of the
parameter the Robin boundary condition begins to behave like a
Dirichlet boundary condition, meaning the potential is less free to
vary which leads to less variance in the inferred parameter in these
regions.

Finally, the centre images in Figures \ref{fig: var_reduc} and
\ref{fig: var_reduc_an} further illustrate the extent to which the
conventional error model leads to overly optimistic (narrow)
confidence intervals. This feature is especially evident in the
anisotropic case (Figure \ref{fig: var_reduc_an}), where the posterior
variance using the conventional error model is significantly smaller
than the reference posterior variance.

\begin{figure}[h!]
\includegraphics[width=\linewidth]{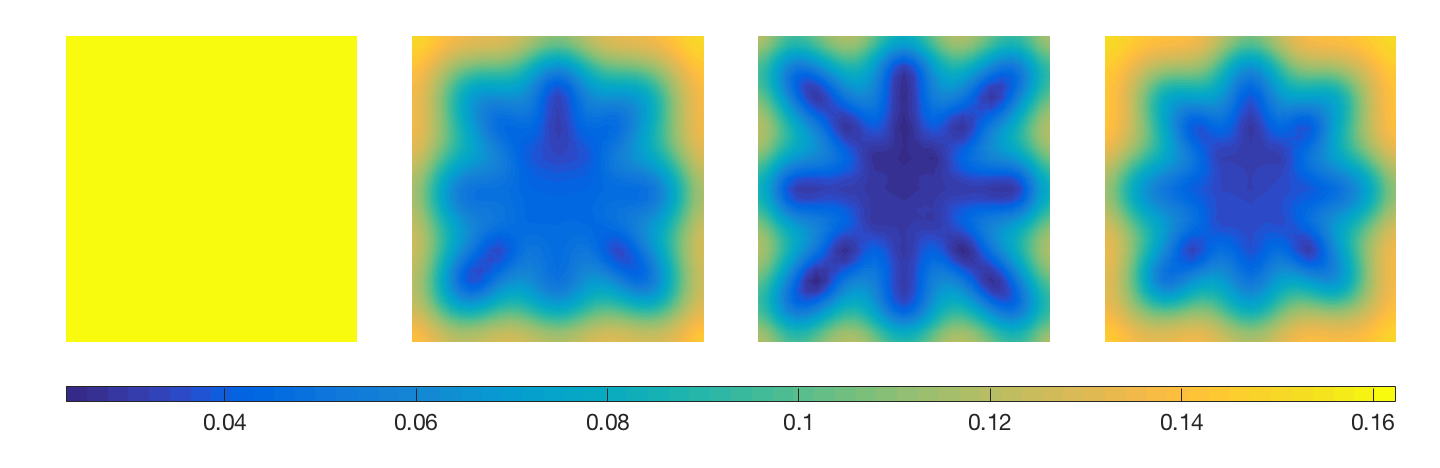}
\caption{The reduction in variance for the anisotropic case. Far left:
  The diagonal of the prior covariance matrix, $\bs{\Gamma}_\beta$, as
  outlined in section \ref{sec: BA2IDIP}. Centre left: The diagonal of
  the reference approximate posterior covariance matrix,
  $\bs{\Gamma}^{\rm REF}_{\rm post}$. Centre right: The diagonal of
  the conventional error model approximate posterior covariance
  matrix, $\bs{\Gamma}^{\rm CEM}_{\rm post}$. Far right: The diagonal
  of the approximation error model approximate posterior covariance
  matrix, $\bs{\Gamma}^{\rm BAE}_{\rm post}$.}\label{fig:
  var_reduc_an}
\end{figure}

\section{Discussion}
In this paper, we considered the problem of inferring the distributed (basal) Robin
coefficient from surface measurements under an unknown random conductivity field.
The forward model at hand was the (anisotropic) Poisson equation with mixed boundary conditions.
To account for model errors that
stem from the uncertainty in the conductivity
coefficient in the underlying PDE, we carry out approximative marginalization
over the conductivity.
In this process, we approximate the related modelling errors and uncertainties as
normal, which is also referred to as the Bayesian approximation error (BAE) approach.

The uncertainty analysis presented here relies on local linearization of
the parameter-to-observable maps at the MAP point estimates, leading to a normal 
(Gaussian) approximation
of the parameter posterior density, which is also referred to as the Laplace approximation. 
We considered two cases of the conductivity field, an isotropically smooth field and an anisotropically
smooth (horizontal strata) one.
The results indicate that fixing the conductivity as an incorrect but otherwise well justified 
(distributed parameter) field can result in infeasible and misleading posterior estimates in the sense 
that the true parameter is not supported by the posterior model.
On the other hand, carrying out approximative marginalization does provide feasible estimates
with both isotropically and anisotropically smooth unknown conductivities.

The computational feasibility (in large-scale) distributed Robin
coefficient problems is provided by the adjoint method and being able
to avoid the simultaneous estimation of the conductivity, which in
contrast to the Robin coefficient, is a random field in the entire
domain.

Future work will be concentrated in two distinct direction. Firstly from a qualitative stand point it would be of interest to consider how the quality of the estimate changes when the number and location of measurements changes. Secondly, in the simplified model (with fixed $a=a_*$) there are no parameters distributed in the domain, warranting an investigation into the computational feasibility of other numerical methods such as boundary element methods (BEM) for solving the the simplified model.

\bibliographystyle{IJ4UQ_Bibliography_Style}

\end{document}